\documentclass[11pt,a4paper,leqno]{amsart}

\usepackage{amsmath,amssymb,amsthm}
\usepackage[dvips]{graphicx}
\usepackage[dvips]{hyperref}

\def\dis{\displaystyle}

\def\R{{\mathbb R}}

\def\virgp{\raise 2pt\hbox{,}}

\def\({\left(}
\def\){\right)}
\def\<{\left\langle}
\def\>{\right\rangle}

\def\d{{\partial}}
\def\e{\varepsilon}
\def\l{\lambda}
\def\si{{\sigma}}

\def\v{{\tt v}}

\theoremstyle{plain}
\newtheorem{theo}{Theorem}[section]
\newtheorem{lem}[theo]{Lemma}
\newtheorem{cor}[theo]{Corollary}
\newtheorem{prop}[theo]{Proposition}

\newtheorem{question}{Question}
\theoremstyle{definition}
\newtheorem{defin}{Definition}[section]

\theoremstyle{remark}
\newtheorem{rema}[theo]{Remark}
\newtheorem*{rema*}{Remark}

\numberwithin{equation}{section}

\begin{document}
\DeclareGraphicsExtensions{.eps} 

\title[Blow-up
  time for NLS]{Monotonicity properties of blow-up
  time for nonlinear Schr\"odinger equation: numerical tests}
\author[C. Besse]{Christophe Besse}
\address[C. Besse]{Simpaf Team -- Inria Futurs\\ Laboratoire Paul
  Painlev\'e, UMR CNRS 8524\\
Universit\'e des Sciences et Technologies de Lille\\
Cit\'e Scientifique\\ 59655 Villeneuve d'Ascq Cedex\\ France}
\email{Christophe.Besse@math.univ-lille1.fr}
\author[R. Carles]{R{\'e}mi Carles}
\address[R. Carles]{MAB, UMR CNRS 5466\\
Universit{\'e} Bordeaux 1\\ 351 cours de la Lib{\'e}ration\\ 33405 Talence
cedex\\ France}
\email{Remi.Carles@math.cnrs.fr}
\author[N. Mauser]{Norbert J. Mauser}
\address[N. J.  Mauser]{Wolfgang Pauli Institute c/o Fak. f. Math. \\
Universit\"at Wien\\
Nordbergstr. 15 \\ A 1090 Wien\\ Austria}
\email{mauser@courant.nyu.edu}
\author[H. P. Stimming]{Hans~Peter Stimming}
\address[H. P. Stimming]{Wolfgang Pauli Institute
and Fakult\"at f. Math. \\ Universit\"at Wien\\ 
Nordbergstr. 15 \\ A 1090 Wien\\ Austria}
\email{hans.peter.stimming@univie.ac.at}
\thanks{
This work was supported by the Austrian Ministry of Science
(BM:BWK) via its grant for the Wolfgang Pauli Institute and by the
Austrian Science Foundation (FWF) via the START Project
(Y-137-TEC) and by the European network HYKE
funded by the EC as contract HPRN-CT-2002-00282.
}
\begin{abstract}
We consider the focusing nonlinear Schr\"odinger equation, in the
$L^2$-critical and supercritical cases. We investigate numerically
the dependence of the blow-up time on a parameter in three cases:
 dependence upon the coupling constant, when the initial data
are fixed;  dependence upon the strength of a quadratic
oscillation in the initial data when the equation and the initial
profile are fixed; finally, dependence upon a damping factor
when the initial data are fixed. It turns out that in most
situations monotonicity in the evolution of the blow-up time
does not occur. In the 
case of quadratic oscillations in the initial data, with
critical nonlinearity, monotonicity holds; this is proven
analytically. 
\end{abstract}

\subjclass[2000]{35Q55, 65M70, 81Q05}
\maketitle


\section{Introduction}
\label{sec:intro}

Consider the Schr\"odinger equation with focusing 
nonlinearity ($\lambda,\si>0$): 
\begin{equation}
  \label{eq:nls}
  i\d_t u +\Delta u = -\lambda |u|^{2\si}u\ , \ (t,x)\in
  \R_+\times \R^n\quad ;\quad
  u_{\mid t =0 } = u_0\, .
\end{equation}
Such an equation may appear as an envelope
equation in the propagation of lasers (see e.g. \cite{Sulem}, and 
\cite{DJMR,JMRIndiana} for a rigorous mathematical justification). It
is well known that if $u_0\in H^1(\R^n)$ (inhomogeneous Sobolev space)
and $\si <\frac{2}{n-2}$, then \eqref{eq:nls} has a unique solution in
$H^1(\R^n)$, defined locally in time (see e.g. \cite{Caz}). It
needs not remain in $H^1(\R^n)$ globally in time: finite time blow-up
may occur when $\si \geq \frac{2}{n}$ and $\lambda >0$ (focusing, or
attractive, nonlinearity); see
e.g. \cite{BourgainLivre,Caz,Sulem}. Since $\lambda\in \R$, the
$L^2$-norm of $u(t,\cdot)$ is independent of time, and finite time
blow-up means that there exists 
$T>0$ (we consider only forward time evolution) such that:
\begin{equation*}
  \|\nabla_x u(t)\|_{L^2}\to +\infty \quad \text{as }t\to T\, .
\end{equation*}
Many papers provide important properties about the blow-up rate or the
profile at blow-up; see
e.g.
\cite{Weinstein83,Weinstein86,MerleCMP,MerleDuke,BourgainWang,Nawa99,
  Perelman,MerleRaphaelInvent,MerleRaphaelCMP,RaphaelMA} for some
analytical results, and \cite{PierreXEDP} for a nice survey of the latest
results. We also point out that recent numerical experiments
\cite{FibichGavishWang} have
shown a completely new phenomenon, where the blow-up profile is
independent of the usual ground state (quasi self-similar ring profile), and
the blow-up rate seems to be the minimal one given by the theory
(square root blow-up rate). As recalled in \cite{PierreXEDP}, 
the main three directions of research in this subject are: 
giving sufficient conditions to have finite time blow-up in the
energy space; estimating the blow-up rate and the stability of the
blow-up r\'egimes; describing the spatial structure of the
singularity formation. 
On the other hand, it  
seems little attention has been paid to the time where blow-up
occurs. In this paper, we investigate by numerical experiments  the
dependence of the blow-up time upon, for
instance, the coupling constant $\lambda$, when the initial
datum $u_0$ is fixed. 

\smallbreak

To motivate our study, we recall some results from
\cite{CW92} and \cite{Fibich01}. In \cite{CW92}, the authors prove
that if the initial datum $u_0(x)$ is replaced by $u_0(x)e^{-ib|x|^2/4}$,
then one can relate explicitly the blow-up time of the corresponding
new solution $u_b$ to that of $u$, in the case of a critical
nonlinearity, $\si=\frac{2}{n}$. This is a consequence of the
conformal invariance. In the super-critical case $\si>\frac{2}{n}$, the
conformal transform does not leave \eqref{eq:nls} invariant, and
introduces a factor $(1-bt)^{n\si -2}$ in front of the
nonlinearity. It is also established in \cite{CW92} that if $u$ has
negative energy (in this case, there is finite time blow-up at least
if $xu_0 \in L^2(\R^n)$ \cite{Glassey}; see also
\cite{OgawaTsutsumi,Martel97} for weaker assumptions), then \emph{for
  large $b$}, 
blow-up occurs sooner than for $b=0$; unlike in the conformally
invariant case, one does not know whether the blow-up time is
monotonous with respect to $b$. The numerical experiments we present
here show that it is not monotonous with respect to $b$. 

\smallbreak

In \cite{Fibich01}, the author considers the damped cubic Schr\"odinger
in space dimension two:
\begin{equation}
  \label{eq:dampedcubic}
  i\d_t \psi +\Delta \psi = -i\delta \psi- |\psi|^{2}\psi\
  , \ (t,x)\in 
  \R_+\times \R^2\quad ;\quad
  \psi_{\mid t =0 } = u_0\, .
\end{equation}
It is conjectured that the blow-up time is monotonous with respect to
$\delta>0$, and guessed that the same holds in the super-critical
case. Our numerical experiments show that neither of the two guesses
is satisfied. These two guesses are very satisfactory when one think
of the initial datum $u_0$ as a single hump; numerics also suggest that
in this case, the blow-up time is monotonous with respect to
$\delta$. On the other hand if $u_0$ is made of, say, two humps, then
the initial intuition seems to be wrong.  

Note that introducing $u(t,x)= e^{\delta t}\psi (t,x)$,
Equation~\eqref{eq:dampedcubic} is equivalent to:
\begin{equation}
  \label{eq:dampedcubic2}
  i\d_t u +\Delta u = - e^{-2 \delta t}|u|^{2}u\
  , \ (t,x)\in 
  \R_+\times \R^2\quad ;\quad
  u_{\mid t =0 } = u_0\, .
\end{equation}
As in the case of initial quadratic oscillations mentioned above, this
transform yields an equation of the form \eqref{eq:nls}, with a
time-dependent coupling ``constant'', $\lambda = e^{-2 \delta
  t}$. This function of time is monotonous, decreasing. 

\smallbreak

This naturally
leads us to the question, when $\lambda$ is really a constant: is the
blow-up time monotonous with respect to $\lambda$?
Tracking the dependence of the blow-up time upon $\lambda$ can be
viewed as both a generalization and a study case of the two problems
raised in \cite{CW92} and \cite{Fibich01}. Numerics show that both in
the critical ($\si = \frac{2}{n}$) and in the supercritical case
($\si>\frac{2}{n}$), one should not expect the blow-up time to be
monotonous with respect to $\lambda$. Essentially, the idea is the
same as what was announced above: when the initial datum $u_0$ is a
single Gaussian, then monotonicity seems to holds. When $u_0$ is the
superposition of two such functions, there is a lack of
monotonicity. 

\smallbreak

Roughly speaking, suppose that two humps are 
placed to both sides of the origin and 
that the natural (free) time evolution tends to send
mass from both these humps to the origin.
This may occur thanks to a particular phase term, for instance, or
simply mass dispersion. 
For large $\lambda$, the blow-up occurs before
the two humps have merged: the two humps do not
interact before the blow-up, and break down independently. For small
$\lambda$, these 
humps merge into one hump before blow-up takes place.
Asymptotically,
for ``large'' $\lambda$ or for ``small'' $\lambda$, we observe some
monotonicity in the blow-up time. On the other hand, the monotonicity
breaks down in the ``transition'' region, that is, for intermediate
$\lambda$. 

\smallbreak

The rest of this paper is organized as follows. In Section~\ref{sec:theory},
we recall some analytical results on local existence and finite time
blow-up for the nonlinear Schr\"odinger equations that appear in this
paper. For the sake of readability, proofs are given in an
appendix. Numerical tests are presented in Section~\ref{sec:numerconst}
to \ref{sec:numerdamped}. In Section~\ref{sec:numerconst}, we consider
the dependence of the blow-up time upon $\lambda$ in \eqref{eq:nls};
in Section~\ref{sec:numerquadr}, we measure the dependence of the
blow-up time upon a modification of the initial datum with quadratic
oscillations; Section~\ref{sec:numerdamped} is devoted to tests on
\eqref{eq:dampedcubic2} and its natural generalization in one space
dimension.


\section{Some theoretical results}
\label{sec:theory}

In this section, we present some analytical results that provide
bounds, from above and/or from below, for blow-up time. The techniques
are classical: the proofs rely on Strichartz estimates,
conservations of mass and energy, and pseudo-conformal conservation
law. Yet, it seems that the explicit 
dependence of the existence time upon some parameters had not been
investigated before, except in \cite{CW92} (see
Section~\ref{sec:quadgene}). These results are somehow a more
quantitative motivation for the numerical tests that follow. 
Technical proofs are given in Appendix~\ref{sec:proofs}. 
\subsection{Generalities}
\label{sec:general}

We first state a local existence result, from which we infer a lower bound on
the blow-up time. Recall that in \eqref{eq:nls}, we consider only
nonnegative time. 
\begin{prop} \label{prop:local}
  Let $\lambda> 0$, $\si \geq \frac{2}{n}$ with $\si < \frac{2}{n-2}$ if $n\geq
3$, and $u_0\in H^1(\R^n)$. There exists $\e >0$ independent of
$\lambda, \si$ and $u_0$ such that if 
\begin{equation}
  \label{eq:cond}
  \lambda T^{\frac{2-(n-2)\si}{2\si +2}} \|u_0\|_{L^2}^{\frac{\si}{\si
    +1}( 2-(n-2)\si)}
\|\nabla_x u_0\|_{L^2}^{\frac{n\si^2}{\si
    +1}}\leq \e\, ,
\end{equation}
then   \eqref{eq:nls}
has a unique solution $u\in C\([0,T[;H^1\)\cap
L^{\frac{4\si+4}{n\si}} \([0,T[; L^{2\si+2}\)$.
If moreover 
\begin{equation*}
  u_0 \in \Sigma := \left\{ f\in H^1(\R^n)\quad ;\quad |x|f\in
  L^2(\R^n)\right\}\, ,
\end{equation*}
then $ u \in C\([0,T[;\Sigma\)$. In addition, the following
quantities are independent of time:
\begin{align}
    &\textrm{Mass: } M=\|u(t)\|_{L^2} ={\rm Const}=
    \|u_0\|_{L^2},\label{eq:mass} \\  
    &\textrm{Energy: } E=\|\nabla_x u(t)\|_{L^2}^2-
\frac{\lambda}{\si +1} \|u(t)\|_{L^{2\si +2}}^{2\si +2} ={\rm
    Const}.\label{eq:energy} 
  \end{align}
\end{prop}
This result is standard (see e.g. \cite{Caz}), except that usually
Condition~\eqref{eq:cond} is not given explicitly: see
Section~\ref{sec:proplocal} for the proof. 
Note that in 
Corollary~\ref{cor:improved} below, we modify this condition, to a
somehow weaker one. The reason is that in
Proposition~\ref{prop:local}, local solutions are constructed without
taking the conservations of mass and energy into account. We shall not
recall explicitly the fact 
that for small initial data, the solution to \eqref{eq:nls} does not
blow up. Indeed, for small data, the conservations of mass and energy
yield an {\it a priori} bound on the $H^1$-norm of the solution, thus
ruling out finite time blow-up. In our case, small initial data means
that, for instance, 
$\lambda^{1/(2\si)}\|u_0\|_{H^1}\leq \delta$ for some constant $\delta$
depending only on $n$ and $\si$. We shall not
insist on that aspect, since our analysis is focused on regimes
where blow-up does occur. This is why in the next two corollaries,
$\l$ is morally ``large'', while  $u_0$ is fixed. 
\begin{cor}\label{cor:improved}
 Let $\lambda> 0$, $\si \geq \frac{2}{n}$ with $\si < \frac{2}{n-2}$ if $n\geq
3$, and $u_0\in H^1(\R^n)$. There exists $\e >0$ independent of
$\lambda$ such that if 
\begin{equation}
  \label{eq:cond2}
  \(\lambda^{2\si}  
\|\nabla_x u_0\|_{L^2}^{4\si}+1\)T^{2-(n-2)\si}\leq \e\, ,
\end{equation}
then  \eqref{eq:nls} has a unique solution $u\in C\([0,T[;H^1\)\cap
L^{\frac{4\si+4}{n\si}} \([0,T[; 
L^{2\si+2}\)$.
If moreover  $u_0 \in \Sigma$,
then $ u \in C\([0,T[;\Sigma\)$.
\end{cor}
The proof of Corollary~\ref{cor:improved} is given in
Section~\ref{sec:corimproved}. 
\begin{cor}[Dependence with respect to the coupling constant]
 Let $\lambda >0$, $\si \geq \frac{2}{n}$ with $\si < \frac{2}{n-2}$ if $n\geq
3$, and $u_0\in \Sigma$. Assume that $u$ blows up in finite time
$T^*>0$.\\
$1.$ We have $\dis T^* \geq C\<\l\>^{-\frac{2\si}{2-(n-2)\si}}$, for
some constant $C$ independent of $\l$, where $\< \l\>=\sqrt{1+\l^2}$.\\
$2.$ If in addition $E<0$, then $\dis T^* \leq C' \<\l\>^{-1/2}$, for
some constant $C'$ independent of $\l$.
\end{cor}
\begin{rema*}
  The above constants $C$ and $C'$ are independent of $\l$, but depend
  on the other parameters,  $u_0$, $n$ and $\si$. As a matter of fact,
  the construction we use yields constants depending only on
  $\|u_0\|_\Sigma$, $n$ and $\si$. 
\end{rema*}
\begin{proof}
  The first part is a straightforward consequence of
  Corollary~\ref{cor:improved}. The second follows from the
  Zakharov--Glassey method \cite{Zakharov,Glassey}. Introduce
  \begin{equation*}
    y(t) = \int_{\R^n}|x|^2 |u(t,x)|^2dx\, .
  \end{equation*}
Then since $\lambda >0$ and $\si \geq \frac{2}{n}$, we have:
$  \ddot y(t) \leq 8n\si E$,
where $E$ denotes the energy defined in
Proposition~\ref{prop:local}. Integrating, we infer that for large
$\lambda$, $y(t) \leq -C\l t^2$, for some positive $C$ independent of
$\lambda$. Since $y(t)\geq 0$ so long as $u$ 
remains in $\Sigma$, this yields the second point of the corollary.
\end{proof}

  We always have $\frac{2\si }{2-(n-2)\si}  > \frac{1}{2}$, so the
  above two bounds go to zero with different rates when  $\l\to
  +\infty$. Condition~\eqref{eq:cond} would yield only $\dis
  T^* \geq C\<\l\>^{-\frac{2\si+2}{2-(n-2)\si}}$. Without even
  trying to see if any of these bounds is sharp, we ask the
  following question:
  \begin{question}
    For $\si \geq \frac{2}{n}$ and a fixed
     initial datum 
$u_0\in \Sigma$,  is the blow-up time for $u$  solution to
\eqref{eq:nls} monotonous with respect to $\lambda$? 
  \end{question}
  
This issue is addressed numerically in
Section~\ref{sec:numerconst}, where our results show that the answer
to the above question should be no.
\subsection{Initial data with quadratic oscillations}
\label{sec:quadgene}
Like in \cite{CW92}, we now fix the equation, and alter only the
initial data, with quadratic oscillations:
\begin{equation}
  \label{eq:nlsCW0}
  i\d_t u + \Delta u = - |u|^{2\si} u \, ,\ x\in
  \R^n \quad ; \quad
  u_{\mid t=0} = u_0(x)\, ,
\end{equation}
with $u_0\in \Sigma$. For $a\not =0$, define:
\begin{equation}\label{eq:chgt}
v(t,x)= \frac{e^{i\frac{|x|^2}{4(t-a)}}}{h(t)^{n/2}}u\left(
\frac{at}{a-t}\virgp \frac{x}{h(t)}\right) , \text{ where }h(t)=
\frac{a-t}{a}\, \cdot
\end{equation} 
Then $v$ solves:
\begin{equation}
  \label{eq:nlsCW}
  i\d_t v + \Delta v = - h(t)^{n\si -2}|v|^{2\si} v  \quad
  ; \quad v_{\mid t=0} = u_0(x)e^{-i\frac{|x|^2}{4a}}\, .
\end{equation}
In the conformally invariant case $\si =2/n$, $v$ solves the same
equation as $u$. The only difference is the presence of (additional)
quadratic oscillations in the data. 
\begin{prop}\label{prop:quadth}
Let $u_0\in\Sigma$ and $2/n\leq \si<2/(n-2)$. Suppose that $u$
blows up at time $T>0$. Let $a\in \R^*$. 
\begin{itemize}
\item If  $a>0$, then $v$ blows up at $T_a(v) = \dis\frac{a}{a+T}T<T$. 
\item If $a< 0$ and $a+T<0$, then $v$ blows up at $T_a(v) =
  \dis\frac{a}{a+T}T>T$. 
\item If $a< 0$ and $a+T\geq 0$, then $v$ is globally defined in
  $\Sigma$ for positive times 
  (but blows up in the past if $a+T> 0$).  
\end{itemize}
\end{prop}
For the critical case $\si =2/n$, this result is proved in \cite{CW92}
(see also \cite{Caz}). We sketch a slightly different proof in
Appendix~\ref{sec:propquadth}, which easily includes the case $2/n<
\si<2/(n-2)$.

In the super-critical case, a natural question is to understand the
role of the function $h$. Introduce $w$ solving:
 \begin{equation}
  \label{eq:nlsCW2}
  i\d_t w + \Delta w = - |w|^{2\si} w  \quad
  ; \quad w_{\mid t=0} = u_0(x)e^{-i\frac{|x|^2}{4a}}\, .
\end{equation}
\begin{prop}
  Let $u_0\in\Sigma$ and $2/n< \si<2/(n-2)$. 
  \begin{itemize}
  \item Let $T_a(w)$ denote the maximal existence time in the future for
  $w$. Then there exists $C=C(\|u_0\|_\Sigma,n,\si)$ independent of
  $a$ such that: 
$$T_a(w) \geq C |a|^{\frac{4\si}{2-(n-2)\si}}.$$
  \item If the energy $E$ of $u$ is negative, then for $a>0$, $w$ blows up at
  time $T_a(w) \leq a$. 
  \end{itemize}
\end{prop}
The first point is a direct consequence of
Corollary~\ref{cor:improved}. The second point is proven in
\cite{CW92}, and relies on the pseudo-conformal law  for
$w$.

To understand the influence of the quadratic oscillations on the
blow-up time, we have to compare the blow-up time of $u$ and that of
$w$. 
In the critical case, the blow-up time depends explicitly on the
magnitude of the quadratic oscillations \emph{via}
Proposition~\ref{prop:quadth}, since $v\equiv w$ by conformal
invariance.  In the super-critical case, we ask:
\begin{question}
  For $\si >\frac{2}{n}$ and a fixed $u_0 \in \Sigma$, is the blow-up
  time for $w$ solving \eqref{eq:nlsCW2} monotonous with respect to $a$?
\end{question}

This issue is addressed numerically in Section~\ref{sec:numerquadr}:
we first compare the numerics with the analytical results in the
conformally invariant case, then perform tests in the supercritical
case
which indicate that the answer to the question above should be no.
\subsection{Damped equation}
\label{sec:dampgene}
 We now consider:
\begin{equation}
  \label{eq:damped0}
  i\d_t \psi + \Delta \psi = - |\psi|^{2\si} \psi-i\delta
  \psi \, ,\quad (t,x)\in 
  \R_+\times\R^n \quad ; \quad
  \psi_{\mid t=0} = u_0\, ,
\end{equation}
for $\delta >0$. A direct application of (the proof of)
Proposition~\ref{prop:local} shows that $\psi$ is defined on $[0,T[$ for
$T \geq C |\delta|^{-1}$. Note that the sign of $\delta$ is irrelevant
at this stage. To take damping effects into account, introduce 
\begin{equation*}
  u(t,x) = e^{\delta t}\psi(t,x)\, .
\end{equation*}
Then $u$ solves:
\begin{equation}
  \label{eq:damped}
  i\d_t u + \Delta u = - e^{-2\si\delta t}|u|^{2\si} u \quad ;
  \quad u_{\mid t=0} = u_0\, .
\end{equation}
For $\delta$ sufficiently large, $u$ is defined globally in time, in
the future:
\begin{prop}\label{prop:damped}
  Let $\si \geq \frac{2}{n}$ with $\si < \frac{2}{n-2}$ if $n\geq
3$, and $u_0\in H^1(\R^n)$. There exists $C=C(\si,n)$ depending only
on $\si$ and $n$ such that if
\begin{equation}
  \label{eq:conddamp}
  \delta^{\si +1} \geq C \| u_0 \|_{L^2}^{\si(2-(n-2)\si)}\| \nabla u_0
  \|_{L^2}^{n\si^2} ,
\end{equation}
then \eqref{eq:damped} has a unique solution $u\in C\(\R_+;H^1\)\cap
L^{\frac{4\si+4}{n\si}} \(\R_+; L^{2\si+2}\)$.
\end{prop}
The proof is given in Section~\ref{sec:propdamped}. 
The following question is addressed numerically in
Section~\ref{sec:numerdamped}:
\begin{question}
  For $\si \geq \frac{2}{n}$ and a fixed $u_0 \in \Sigma$, is the
  blow-up time for $u$ solution to \eqref{eq:damped} monotonous with
  respect to $\delta >0$?
\end{question}
The simulations show that the answer is no.


\section{Numerical Test, dependence on $\lambda$}
\label{sec:numerconst}

We perform numerical tests by using a direct discretization method
for equation \eqref{eq:nls}, respectively \eqref{eq:damped}. 
By this approach no restrictions are imposed on simulations
on the closeness to eventual  blow-up  points in space and time and by
the evolution in places ``not close'' to the blow-up.
We employ two different numerical methods: the
Time-Splitting Spectral method (TSSP), and the Relaxation method (RS).

\smallskip

The TSSP is based on an operator splitting method, the split-step method.
The split-step method is based on a decomposition of the flow of 
the nonlinear equation \eqref{eq:nls} (or \eqref{eq:damped}).
Define the flow $X^t$ as the flow of the linear Schr\"odinger equation
\begin{equation*}
\left \{ 
  \begin{array}{ll}
i \partial_t v +\Delta v=0, & x \in \R^n, \ t>0, \\
v(x,0)=v_0(x), & x \in \R^2,
  \end{array}
\right .
\end{equation*}
and $Y^t$ as the flow of the nonlinear differential equation
\begin{equation*}
\left \{ 
  \begin{array}{ll}
i \partial_t w = -g(t) |w|^{2\sigma} w, & x \in \R^n, \ t>0, \\
w(x,0)=w_0(x), & x \in \R^2.
  \end{array}
\right .
\end{equation*}
where $g(t)=\lambda$ for the case of  \eqref{eq:nls}
and $g(t)=e^{-2\sigma\delta t}$ for \eqref{eq:damped}.
Then the split-step method consists of approximating the exact $u(x,t)$
at each time step by combining the two flows $X^t$ and $Y^t$.
We employ here the 
Strang formula $Z^t_S = X^{t/2}Y^tX^{t/2}$, which is of second order. 
Higher order splitting is also possible.
For a general nonlinear Schr\"odinger equation, a convergence proof of the
splitting method was 
done by Besse \textsl{et al.} in \cite{BBD}. 
For the TSSP, a spectral method is employed to compute the flow $X^t$ of
the free Schr\"odinger equation. 
The flow $Y^t$, which is the flow  of a nonlinear ODE, 
can be computed exactly, since it leaves $|w|^2$ invariant. So its 
integration is straightforward.
The TSSP has proved to be an efficient and
reliable method for NLS type equations. See for example
\cite{BJM,BMS1} for a study of the NLS
 in the semi-classical limit case, and \cite{PeLi,BMS04} 
for a more general numerical study.

For the 2-d calculations, a parallel version of the TSSP scheme
is used on the parallel cluster machine ``Schr\"odinger III''
at the University of  Vienna.

\smallbreak

The Relaxation method (RS) is a discretization of finite difference type 
 \cite{B04}.
It is based on central-difference approximation shifted by a half time-step. 
To formulate the scheme, we rewrite the NLS equation as a system:
\begin{equation*}
\left \{ 
  \begin{array}{rcc}
i \partial_t u  +\Delta u & = & -g(t)\; \psi \; u, \\ 
\psi & = & |u|^{2\sigma} 
  \end{array}
\right .
\end{equation*}
where $g(t)$ is as above. Then the first equation is discretized by
central difference quotient at
the time $t_{n+1/2} = (n + \frac 1 2) \Delta t$, and the second by
central time average at the time $t_n = n \Delta t$.
Let $u^n$ be the approximation at $t=t_n$, then the scheme is
given by
\begin{eqnarray*}
i \frac {u^{n+1}-u^n}{\Delta t}  +\Delta_D 
\left( \frac {u^{n+1}+u^n}{2} \right)
& = & -g(t) \; \psi^{n+1/2} \left( \frac {u^{n+1}+u^n}{2} \right), \\ 
 \frac {\psi^{n+1/2}+\psi^{n-1/2}}{2} & = & |u^n|^{2\sigma} 
\end{eqnarray*}
where $\Delta_D $ denotes a finite difference Laplacian.
It requires only explicit evaluations of the nonlinear term.
It also conserves energy (\cite{B04}).

\smallbreak

To determine whether blow-up is occurring or not,
we calculate the two terms in the energy \eqref{eq:energy}, 
kinetic and potential energy, and look for an increase
of at least four orders of magnitude in both of them. The first time at
which this  
is occurring is assumed to be the blow-up time.
The space and time resolution are chosen sufficiently fine  such that this 
increase is realized.

\subsection{Critical power}
\smallskip

\begin{figure}
  \center
  \includegraphics[width=0.5\textwidth]{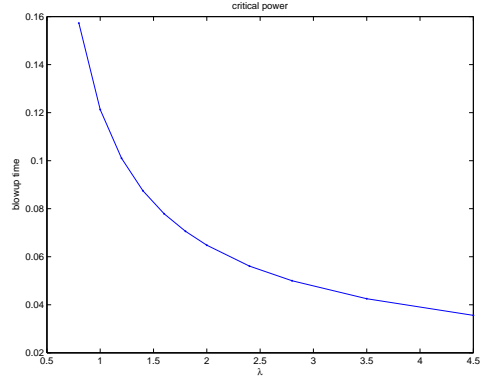}
  \caption{\label{im:test1}
  Blow-up time with varying $\lambda$, single Gaussian data (Test 1).
 }
\end{figure}

\begin{figure}
  \center
  \includegraphics[width=0.7\textwidth]{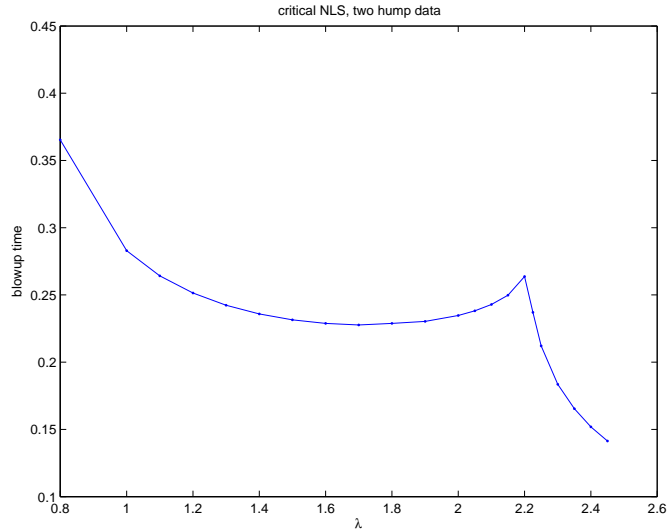}
  \caption{\label{IMAGE01}
  Blow-up time with varying $\lambda$ (Test~2).
 }
\end{figure}

\begin{figure}
  \center
  \includegraphics[width=0.65\textwidth]{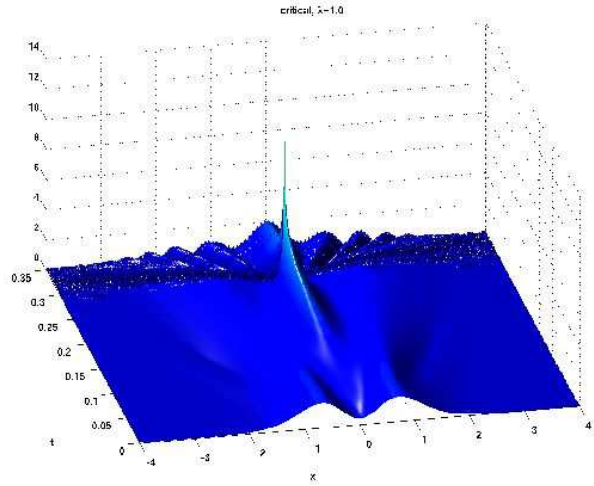}
  \caption{\label{im:lowevo}
  Time evolution with low potential energy.
 }
\end{figure}

\begin{figure}
  \center
  \includegraphics[width=0.65\textwidth]{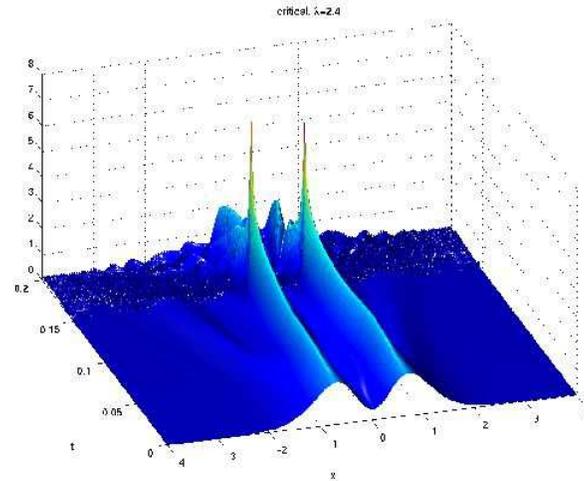}
  \caption{\label{im:highevo}
  Time evolution with high
 potential energy.
 }
\end{figure}

\begin{figure}
  \center
  \includegraphics[width=0.7\textwidth]{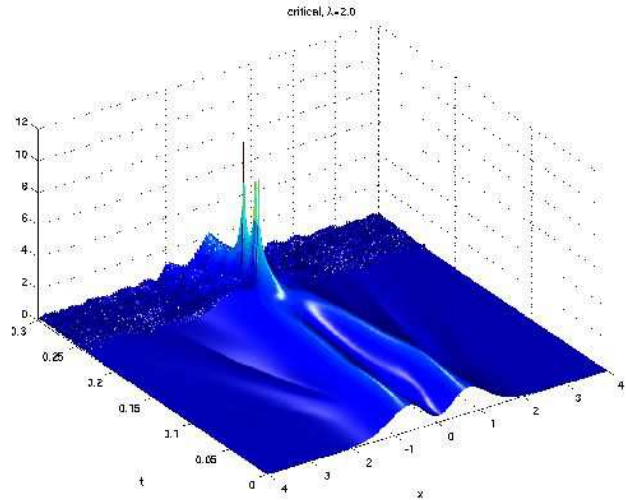}
  \caption{\label{im:medevo}
  Time evolution in intermediate regime: non-monotonicity.
 }
\end{figure}

\begin{figure}[t]
  \center
  \includegraphics[width=0.7\textwidth]{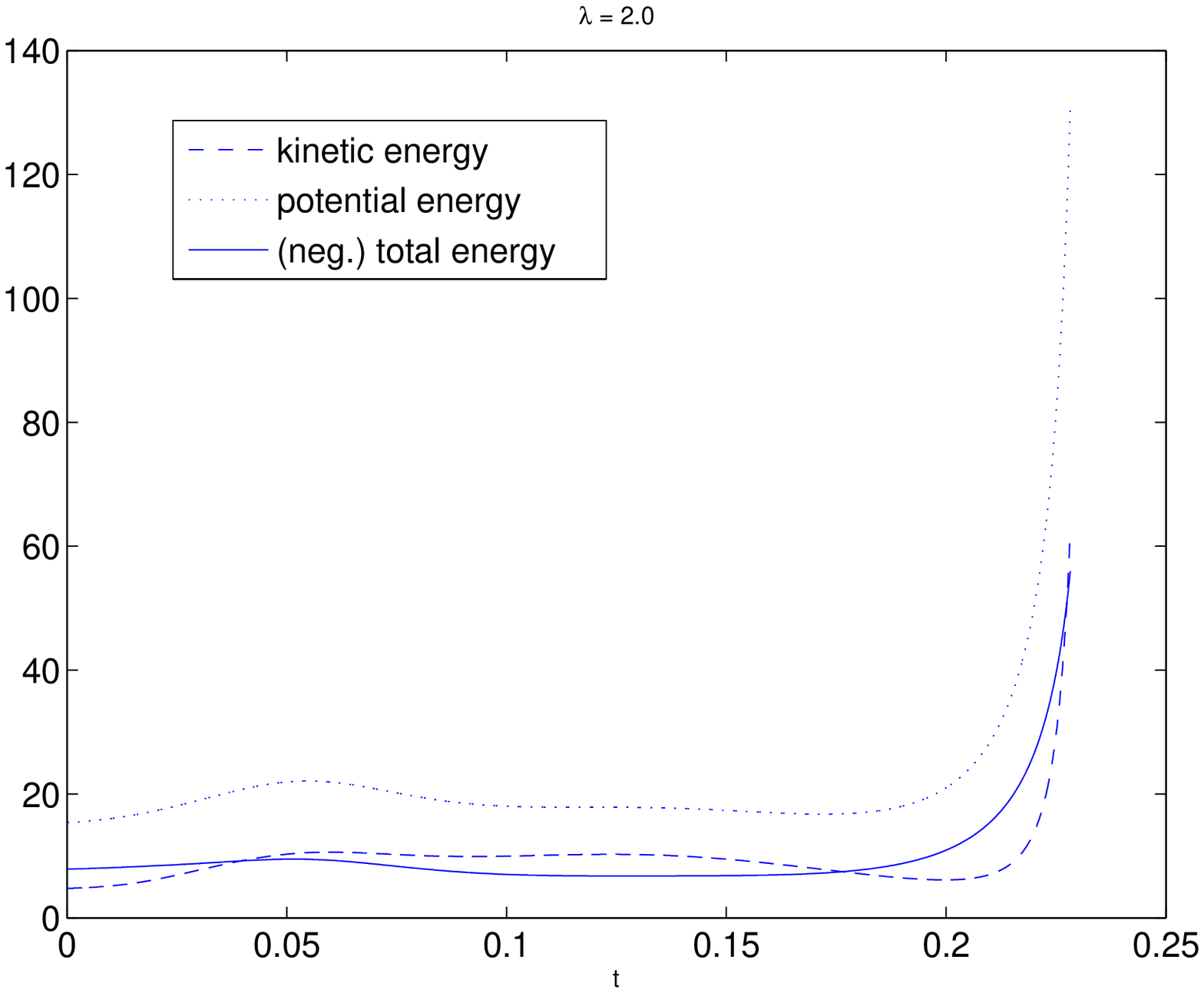} 
  \caption{\label{im:discr_by_enerpt}
  Determining blow-up time by increase of energy: 
 (minus) total energy in comparison with kinetic and
 potential energy away from the blow-up
 }
\end{figure}

First we consider \eqref{eq:nls} for $\sigma = \frac{2}{n}$, 
that is the critical case.
We study the dependence of the blow-up time on the constant $\lambda$ for
a series of different data.

\smallskip
\noindent{\bf Tests in one space dimension}
\smallskip

{\bf Test 1.} For the case $n=1$, the first kind of data we study is
\begin{equation}
\label{eq:singledata}
u_0(x) = C \, e^{-x^2} \, e^{-i\log(e^x+e^{-x})}.
\end{equation}
The constant $C$ is equal to $1.75$ which leads to $\|u_0\|^2_{L^2} = 3.839$.
We take $Np=2^{12}=4096$ mesh points, and several time steps
according to the necessary resolution for
the blow-up, with
$\Delta t  = 2.5 \cdot 10^{-6}$ as smallest.
The discretization domain is $[-8,8]$. 
Figure \ref{im:test1} shows the blow-up time in relation to 
a changing $\lambda$. It can be observed that the blow-up time
is decreasing monotonously with $\lambda$, as predicted by the
heuristics for the case of a single Gaussian profile.
\smallbreak

{\bf Test 2.} The next kind of data we study is
\begin{equation} \label{eq:twohpdata}
u_0(x)=C \, \(e^{-x^2} - 0.9 e^{-3x^2}\) e^{-i\log(e^x+e^{-x})}.
\end{equation}
The constant $C$ is equal to $4$ which leads to $\|u_0\|^2_{L^2} = 3.907$.
The difference of two Gaussian profiles results in two local maxima in the 
modulus of $u_0$. The phase term has a focusing effect and  
its focus point does not agree with the local maxima of the
modulus.
The finest discretization parameters used are
$Np=2^{14}=16384$, $\Delta t = 2.5 \cdot 10^{-6}$.
The discretization domain is $[-8,8]$.

Figure~\ref{IMAGE01} shows the blow-up time in relation to 
a changing $\lambda$. It can be observed that for 
low and very high strengths of the nonlinearity  $\lambda$ 
($\lambda < 1.6$ and  $\lambda > 2.2$), the blow-up time is monotonously 
decreasing with $\lambda$, while in between there is a 
region where monotonicity does not hold.
\smallbreak

Heuristically, two effects  play a role in the solution with this
data: the nonlinear self-focusing, which tries to focus the 
mass to points where the most mass is already present,
and the (linear) phase influence, which tends
to focus the mass at zero.
\smallbreak

In Figures~\ref{im:lowevo}, \ref{im:highevo} and \ref{im:medevo},
the time evolution of the modulus of $u(t,x)$ is shown for 
values of $\lambda$ from the three different regions of the above
curve.\\
$\bullet$ Figure~\ref{im:lowevo} shows the case
$\lambda=1$. The two initial humps  
merge to one hump before the blow-up, which happens at a single
point. The phase focusing happens at a faster time scale than the
nonlinear 
focusing.\\
$\bullet$ Figure~\ref{im:highevo}  shows the case
$\lambda=2.4$. Blow-up  
is occurring simultaneously at two points and there are always two humps 
present. Apparently the nonlinear self-focusing here happens faster than 
phase focusing.\\
$\bullet$ Figure~\ref{im:medevo} is for $\lambda=2.0$, which is in the
non-monotonicity 
region. Blow-up here occurs at a single point, and the merging of the 
two humps is closer to the blow-up than in Figure~\ref{im:lowevo}. 
In this case it is not clear which of the two effects 
would happen at a faster time scale, nor how they would interact.
Blow-up is occurring, but the blow-up time is no longer monotonous 
with respect to the size of the nonlinear term.
\smallbreak

In Figure~\ref{im:discr_by_enerpt}, we show the energy components for the
case $\lambda=2.0$ which are used to determine the blow-up time.
The kinetic energy increases
by a factor of $10^4$, or slightly more, at the blow-up time.
Figure~\ref{im:discr_by_enerpt} 
compares the kinetic and potential energy parts
with $-E(u(t))$ away from the blow-up. The TSSP scheme used for this simulation
does not conserve energy. Near the blow-up time, the energy conservation is no 
longer true for the numerical result.
Note that, outside of a time interval close to blow-up,
the kinetic energy is not monotonously increasing in time.
\smallbreak

\begin{rema*}
If the data \eqref{eq:twohpdata} are used without the phase term,
which leaves just a ``two-hump'' profile, non-monotonicity can be 
observed in the same way as described above.
However overall blow-up times increase. 
Apparently the merging of the two humps can occur even without the  
influence of an initial phase term thanks to the
mass dispersion tendency of the free evolution together with the
focusing effect of the nonlinearity.
\end{rema*}

\begin{figure}
  \center
  \includegraphics[width=0.6\textwidth]{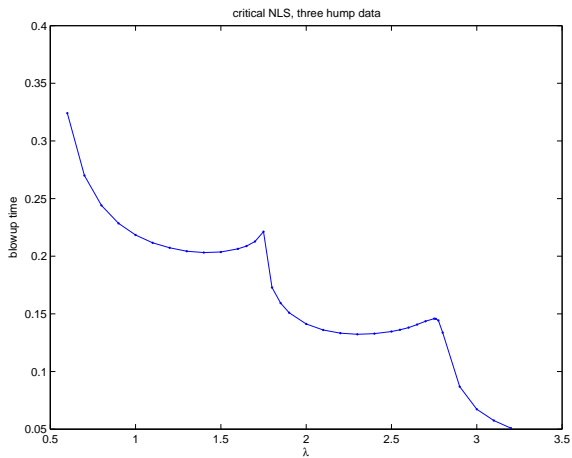}
  \caption{\label{im:threehump}
  Blow-up time with varying $\lambda$, three hump data (Test~3).
 }
\end{figure}

{\bf Test 3.} Data with three humps\\
The next test uses a sum of three Gaussians and the same phase term as before,
so there are three local modulus maxima instead of two.
\begin{equation*}
u_0(x)=C \, \( e^{-(3x)^2} + e^{-(3(x-1))^2} 
+ e^{-(3(x+1))^2}\) e^{-i\log(e^x+e^{-x})},
\end{equation*}
with $C=2$, $\|u_0\|^2_{L^2} = 5.09$.
The discretization parameters are 
$Np=2^{12}$ mesh points and  $\Delta t = 1.5 \cdot 10^{-5}$.
The discretization domain is $[-8,8]$.
The blow-up times with respect to changing $\lambda$ 
are shown in Figure~\ref{im:threehump}.
The same effect as above can be observed.
\smallbreak

There are two regimes for
the nonlinearity strength $\lambda$ where the blow-up time
is non-monotonous. Blow-up happens either at three points, two or one point
depending on $\lambda$.

\smallbreak

\begin{figure}
  \center
  \includegraphics[width=0.6\textwidth]{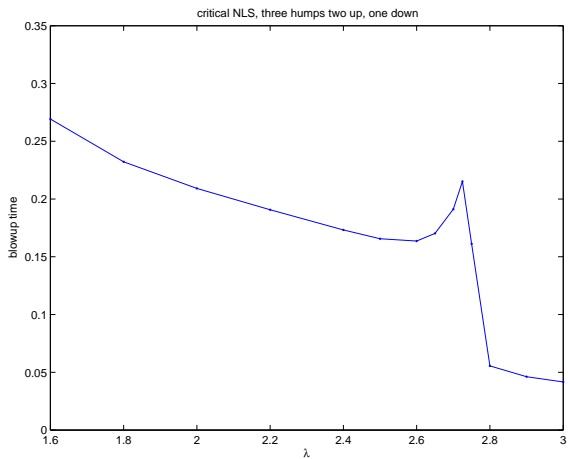}
  \caption{\label{im:tuponedown}
  Blow-up time with varying $\lambda$, three hump data (Test~4).
 }
\end{figure}

{\bf Test 4.} Data with two humps up, one down\\
In this test the data  are taken to be
\begin{equation*}
u_0(x)=C \, \(  e^{-(3(x-1))^2} - e^{-(3x)^2} 
+ e^{-(3(x+1))^2}\) e^{-i\log(e^x+e^{-x})},
\end{equation*}
with $C=2$, $\|u_0\|^2_{L^2} = 4.93$.
One of the three Gaussians has an opposite sign, so there is a constant phase
shift in part of the data. 
The blow-up times are shown in Figure~\ref{im:tuponedown}.
\smallbreak

For $\lambda < 1.6$, there is no blow-up occurring. For larger $\lambda$,
non-monotonicity similar to the situation above can be observed. Observe
that the slope of the curve is rather steep, for
$2.7<\lambda<2.8$. This shows a  highly nonlinear
phenomenon. Note however that there is one point on the steep line after 
the local maximum: 
\begin{center}
\begin{tabular}[c]{c|c}
 $\l$ & $T^*$\\
\hline
$2.7$  &  $0.1912$\\
$2.725$ &   $0.2152$\\
$2.75$  &  $0.1612$ \\
$2.8$  &  $0.0555$ 
\end{tabular}
\end{center}
\smallbreak

\smallbreak

\begin{figure}
  \center
  \includegraphics[width=0.5\textwidth]{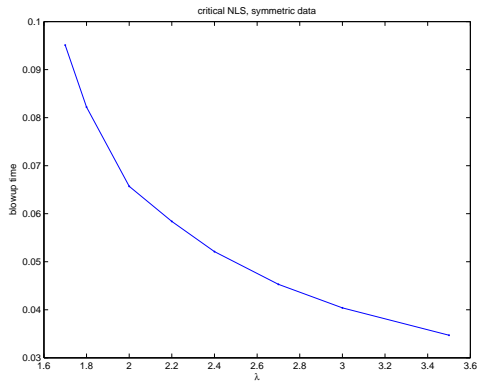}
  \caption{\label{im:pointsymm}
  Blow-up time with varying $\lambda$,  
  data with point symmetry (Test~5).
 }
\end{figure}

{\bf Test 5.} Data with one hump up, one down\\
We use
\begin{equation*}
u_0(x)=C   e^{-x^2} \tanh x  \, e^{-i\log(e^x+e^{-x})}.
\end{equation*}
We use $C=3.0$ and $\|u_0\|^2_{L^2} = 4.4476$, and discretizations of 
$Np=2^{13}$ mesh points and  $\Delta t = 2.0 \cdot 10^{-6}$.
The result is shown in Figure~\ref{im:pointsymm}.
In this case, the blow-up time is monotonous with $\lambda$.

\begin{rema*}
For other point-symmetric data with more than two humps, there
is monotonicity.
\end{rema*}

\smallbreak

{\bf Test 6.} Stability issues: asymmetric data\\
In order to investigate the dependence of the qualitative findings 
on the symmetry of the data, we choose some data which 
perturb the  symmetry of \eqref{eq:twohpdata}. We used
\begin{equation*}
u_0(x)=C \(e^{-x^2} - 0.9 e^{-3x^2}\) e^{-i\log 2 \cosh (x-0.25)},
\end{equation*}
and
\begin{equation*}
u_0(x)=C \(e^{-(x-1.5)^2} + 0.99 e^{-(x+1.5)^2}\) e^{-i\log(e^x+e^{-x})},
\end{equation*}
to have either a phase that focuses away from the symmetry point 
between the two humps, or 
two humps of different heights.
It turns out that non-monotonicity 
can be observed for for both of these cases, too. The perturbations
to \eqref{eq:twohpdata} tested here are small enough not to
hinder  non-monotonous blow-up times.

\bigskip
\noindent {\bf Test in two space dimensions}
\smallskip

\begin{figure}
  \center
  \includegraphics[width=0.7\textwidth]{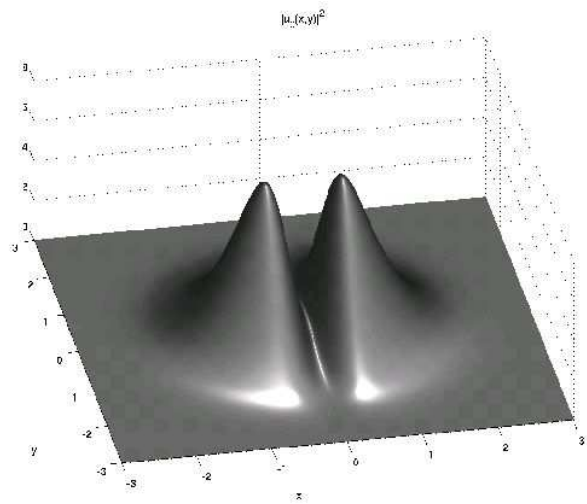}
  \caption{\label{im:tdthdata}
  Initial data for two dimensional case, modulus (Test~7).
 }
\end{figure}

{\bf Test 7.} 
For the test in two space dimensions, we extend \eqref{eq:twohpdata} by
making the phase term radially symmetric and multiplying the
one-dimensional two-hump 
profile by a single Gaussian in the second space dimension:
\begin{equation} \label{eq:tdth_data}
u_0(x,y)=C \(e^{-x^2} - 0.9 e^{-3x^2}\) e^{-y^2} 
e^{-i\log 2 \cosh \(\sqrt{x^2+y^2}\)}.
\end{equation}
Figure \ref{im:tdthdata} shows the modulus of \eqref{eq:tdth_data}.
We choose $C=7.0$, $\|u\|_{L^2}^2=15.0$.
The smallest discretization parameters used are
$Np=2^{12}$ mesh points and  $\Delta t = 1 \cdot 10^{-5}$
with the discretization domain $[-4,4]^2$.

\smallbreak 

\begin{figure}
  \center
  \includegraphics[width=0.6\textwidth]{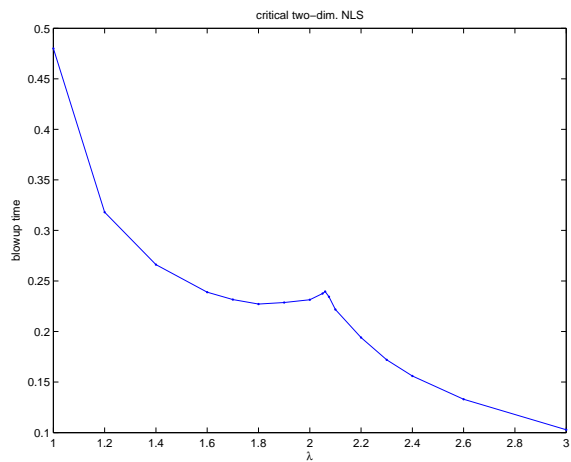}
  \caption{\label{im:twoh}
  Blow-up time with varying $\lambda$, two dimensional case (Test~7).
 }
\end{figure}

The blow-up times with changing $\lambda$ are shown in Figure
\ref{im:twoh}. Non-monotonicity can be observed.

\subsection{Supercritical power}
$\ $

\begin{figure}
  \center
  \includegraphics[width=0.6\textwidth]{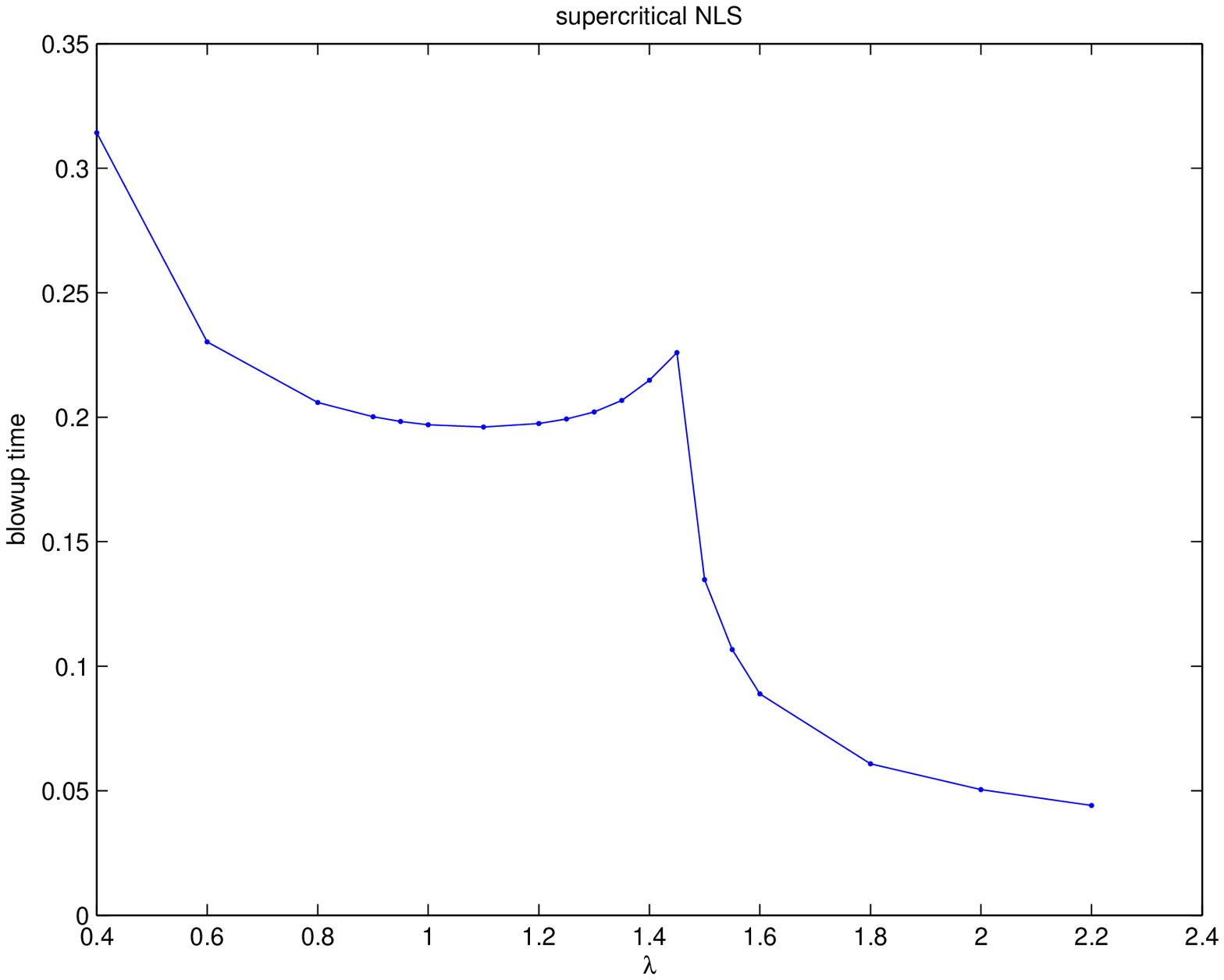}
  \caption{\label{im:supercr}
  Blow-up time with varying  constant $\lambda$ (Test~8).
 }
\end{figure}

{\bf Test 8.}
We tested equation \eqref{eq:nls}  in one space dimension,  with 
$\sigma = 3$ and the data \eqref{eq:twohpdata} ($C=3.5$, hence
$\|u_0\|^2_{L^2} = 2.99$).
The discretization parameters are
$\Delta x = 0.0039$, and up to  $\Delta t = 1 \cdot 10^{-6}$.
The discretization domain is $[-8,8]$.
The blow-up time with varying $\lambda$ is shown in Figure~\ref{im:supercr}.
Also in the supercritical case, non-monotonicity of
blow-up times can be observed.

\section{Numerical Test, dependence on quadratic oscillations}
\label{sec:numerquadr}

We turn to equation  \eqref{eq:nlsCW} and investigate the dependence
of blow-up time on the scale of quadratic oscillations. 
To compare the simulations to the result of Proposition~\ref{prop:quadth},
we simulate \eqref{eq:nlsCW0} with the same data to obtain the 
blow-up time for this equation,
and then 
plot the curve of $T_a$ that is predicted in  Proposition~\ref{prop:quadth}
for the two regions $a > 0$ and $a<0, \ a+T<0$.

\subsection{Critical power}  
$\ $
We use \eqref{eq:nlsCW} in space dimension one,  with $\sigma = 2$
with various $u_0(x)$.

\begin{figure}
  \center
  \includegraphics[width=0.6\textwidth]{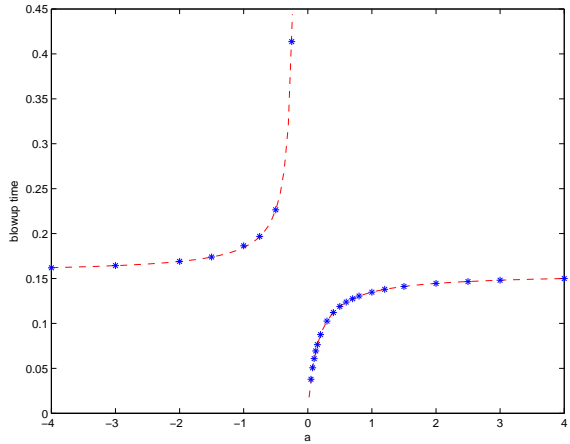}
  \caption{\label{im:quadrsg}
  Blow-up time with varying $a$ in quadratic oscillations,
  critical power (Test~9).
 }
\end{figure}

{\bf Test 9.} Single hump: We take
\begin{equation*}
u_0(x)=C \, e^{-x^2}
\end{equation*}
with $C=1.75$.
The discretization parameters in this and the two following
tests are
$Np=2^{14}$ mesh points and $\Delta t = 4 \cdot 10^{-6}$.
The discretization domain is $[-8,8]$ for all cases except
the two largest negative $a$ in Test 9 and 10. Here the domain
is extended to $[-18,18]$ and the space resolution is $Np=2^{13}$ resp.
$[-40,40]$ and $Np=2^{14}$ for the largest negative $a$ 
of Test~11.
Figure~\ref{im:quadrsg} shows the blow-up time of $v$ in relation 
to the scale $a$ of  quadratic oscillations in the data.
We use both positive $a$ and negative $a$ with $a+T<0$.
Asterisks denote the blow-up times and
the dashed line shows the result of Proposition~\ref{prop:quadth}
with $T$ obtained by a simulation of \eqref{eq:nlsCW0}.
It can be observed that the results agree very well.
In each of the two regions for $a$ that have been used in this
test, the blow-up time  is monotonous with respect to $a$.

\begin{figure}
  \center
  \includegraphics[width=0.6\textwidth]{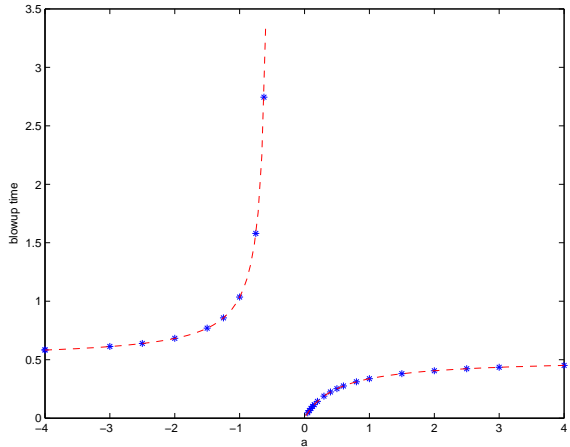}
  \caption{\label{im:quadrdouble}
  Blow-up time with varying $a$ in quadratic oscillations,
  critical power, two hump data (Test~10).
 }
\end{figure}

{\bf Test 10.} Two humps: We take
\begin{equation*} 
u_0(x)=C \, \(e^{-x^2} - 0.9 e^{-3x^2}\)
\end{equation*}
which is the same as \eqref{eq:twohpdata} without the phase
term that appears there. We take $C=4.0$.
Figure~\ref{im:quadrdouble} shows the blow-up times marked by asterisks.
As above the dashed line shows the modified blow-up time according to
Proposition~\ref{prop:quadth}. The two curves agree.

\begin{figure}
  \center
  \includegraphics[width=0.6\textwidth]{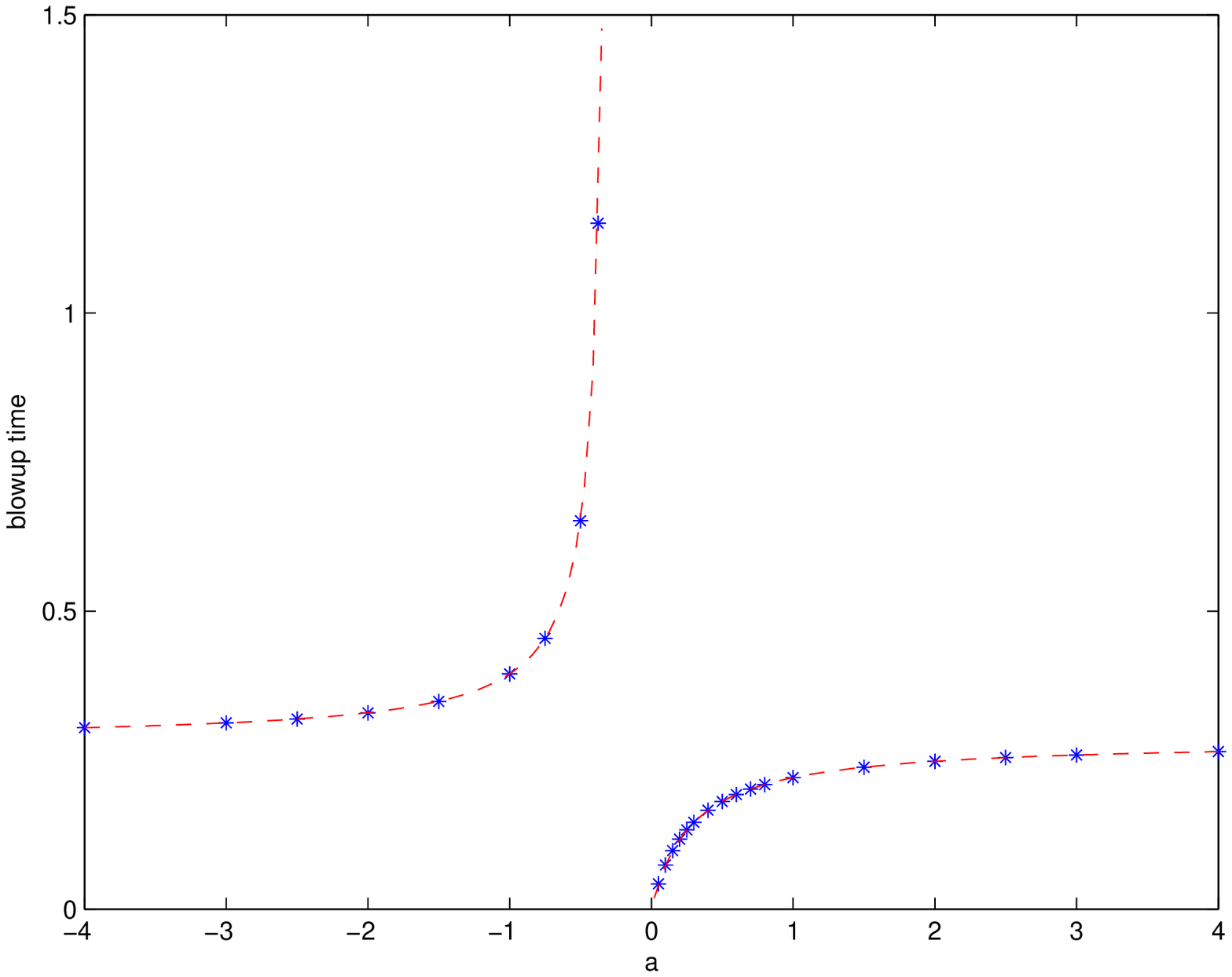}
  \caption{\label{im:quadrphase}
  Blow-up time with varying $a$ in quadratic oscillations,
  critical power, data with additional phase (Test~11).
 }
\end{figure}

{\bf Test 11.} Two humps with additional phase: Here we use  
\begin{equation*}
u_0(x)=C \, \(e^{-x^2} - 0.9 e^{-3x^2}\) e^{-i\log(e^x+e^{-x})}
\end{equation*}
which is the same as \eqref{eq:twohpdata}.
Figure \ref{im:quadrphase} shows the blow-up times of the simulations
and according to the Proposition as above. 
The results agree.


\subsection{Supercritical power}  
$\ $
We now consider Eqns.~\eqref{eq:nlsCW} and \eqref{eq:nlsCW2}
with $\sigma = 3$, in space dimension one. 
We give results for both $w$ and $v$.
We use the same series of data as for the critical power case.

The discretization parameters in Tests 12 to 14  are
$Np=2^{13}$ mesh points and  $\Delta t = 1 \cdot 10^{-5}$, down to
$\Delta t = 4 \cdot 10^{-6}$.

\begin{figure}
  \center
  \includegraphics[width=0.45\textwidth]{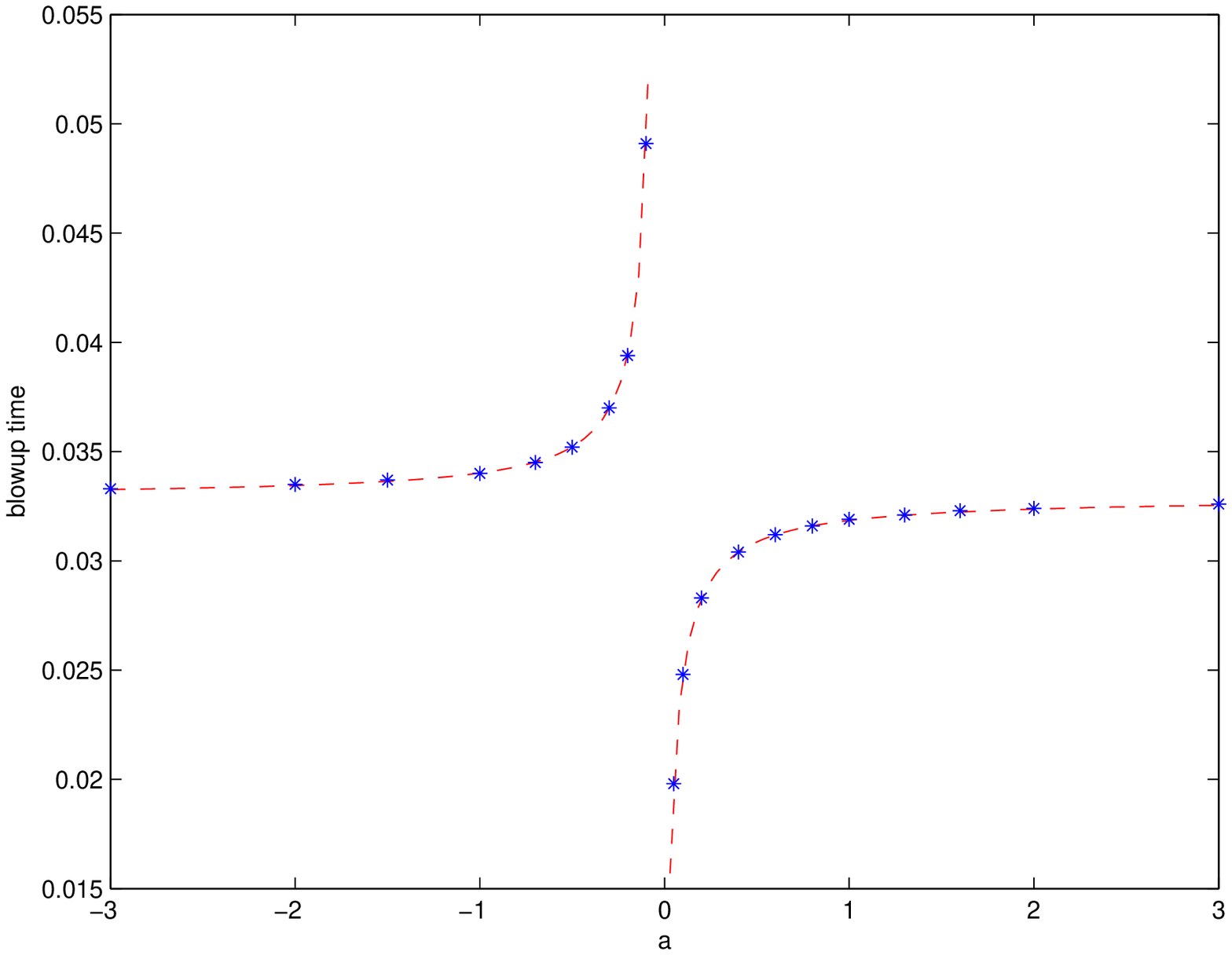}
  \includegraphics[width=0.45\textwidth]{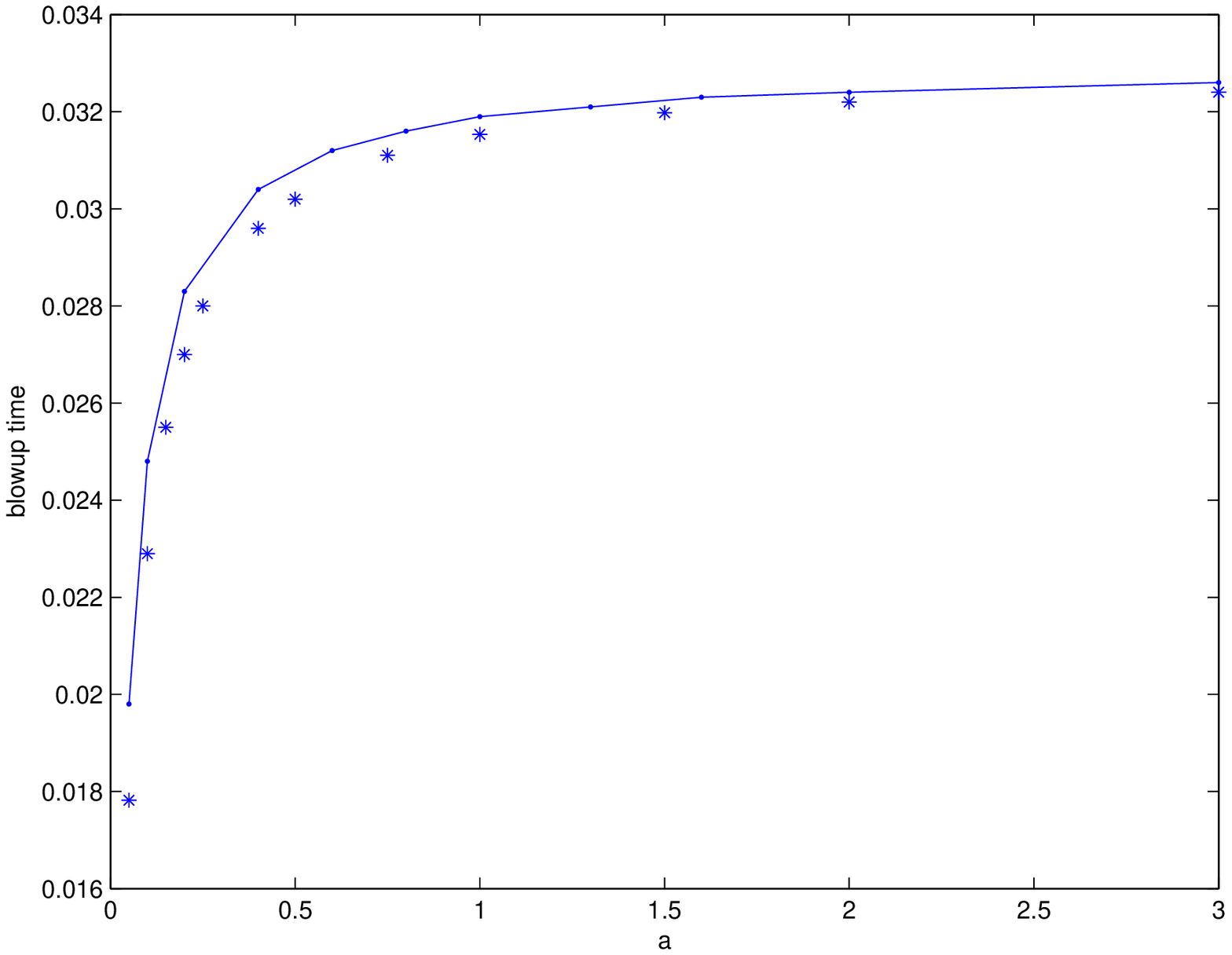}
  \caption{\label{im:quadscr1}
  Blow-up time with varying $a$ in quadratic oscillations,
  supercritical power. Left: blow-up times for $v$.
  Right: comparison of $v$ and $w$, dots with line for $v$,
  asterisks for $w$ (Test~12).
 }
\end{figure}

{\bf Test 12.} Single hump: We take
$
u_0(x)=C \, e^{-x^2}.
$
Figure~\ref{im:quadscr1} shows the blow-up time in relation 
to the scale $a$ of  quadratic oscillations in the data.
In the left figure, the dashed line is the calculated blow-up time for $v$
and the asterisks mark the simulated blow-up times.
The right figure compares the blow-up times for $v$ and $w$, 
where $w$ is denoted by asterisks and $v$ by dots joined by a line.
We can see that $w$ blows up a bit earlier than $v$ for positive $a$.
For negative $a$, it blows up a bit later, but the blow-up times are
rather close 
to those of $v$ as in the positive $a$ case.
Again the blow-up times for $v$ match the result of Proposition 
\ref{prop:quadth}.

\begin{figure}
  \center
  \includegraphics[width=0.5\textwidth]{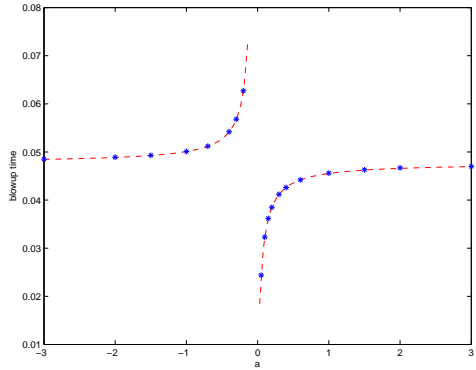}
  \caption{\label{im:quadscr2}
  Blow-up time of $v$ with varying $a$ in quadratic oscillations,
  supercritical power  (Test~13).
 }
\end{figure}

{\bf Test 13.} Two humps: We take
$
u_0(x)=C  \(e^{-x^2} - 0.9 e^{-3x^2}\).
$
Figure~\ref{im:quadscr2} shows the blow-up time of $v$ in relation 
to $a$.

\begin{figure}
  \center
  \includegraphics[width=0.5\textwidth]{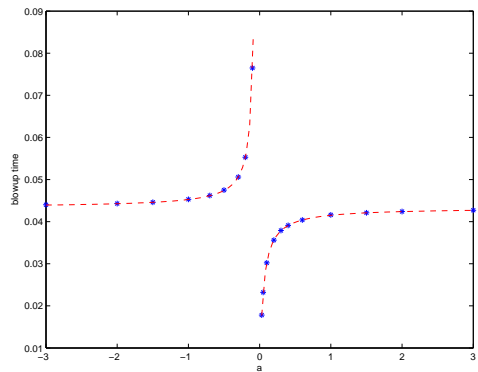}
  \caption{\label{im:quadscr3}
  Blow-up time of $v$ with varying $a$ in quadratic oscillations,
  supercritical power, data with additional phase (Test~14).
 }
\end{figure}

{\bf Test 14.} Two humps with additional phase: Here we use  
\begin{equation*}
u_0(x)=C  \( e^{-x^2} - 0.9 e^{-3x^2} \)  e^{-i\log(e^x+e^{-x})}.
\end{equation*}
Figure \ref{im:quadscr3} shows the blow-up time of $v$ in relation 
to $a$.
\smallskip

\begin{figure}
  \center
  \includegraphics[width=0.7\textwidth]{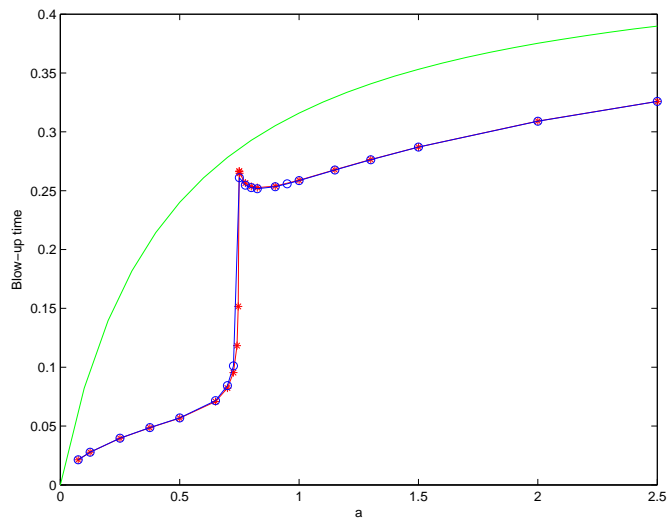}
  \caption{\label{im:quadsc_nom}
  Blow-up time for $w$ with varying $a$ in quadratic oscillations,
  supercritical power, asymmetric data (asterisks and circles). 
  Solid line: Blow-up time for $v$.
 (Test~15).
 }
\end{figure}

{\bf Test 15.} Two humps placed asymmetrically: Here we use  
\begin{equation*}
u_0(x)=C  \( e^{-(3x)^2} + e^{-(3(x-1.5))^2 } \).
\end{equation*}
We take $C=1.8$.
The smallest discretization parameters used in
this test are
$Np=2^{15}$ mesh points and  $\Delta t = 1.5 \cdot 10^{-6}$.
The computation domain is $[-8,8]$, except for the  largest value
of $a$. 
Figure~\ref{im:quadsc_nom} shows the blow-up time for $w$ in relation 
to $a$ obtained by the two numerical methods employed: the 
circles represent simulations done by the TSSP, the asterisks
simulations by the RS. In addition a solid line displays the 
blow-up times for $v$.
Non-monotonicity can be observed, which answers question 2
in a negative way. 
The blow-up time of $w$ is always smaller than that of $v$.
We also see that the results of the two different
schemes agree in a good way. 
Note that the occurrence of non-monotonicity is very sensitive to the
size of the data. If we choose $C = 1.7$ or $1.9$ instead of
$C=1.8$, monotonicity can be observed.

\begin{rema}[Test~15]\label{rema:C}
Leave out the question of quadratic
oscillations, and consider $u$ solving \eqref{eq:nlsCW0}.
Its blow-up time  as
a function of the constant $C$ is given by:
\begin{center}
\begin{tabular}[c]{c|c}
 $C$ & $T^*$\\
\hline
$1.795$ & $0.528$\\
$1.798$ & $0.480$\\
$1.8$ &   $0.462$\\
$1.804$ & $0.446$\\
$1.808$ & $0.507$\\
$1.81$  &  $0.076$ \\
$1.82$  &  $0.048$
\end{tabular}
\end{center}
Since changing $C$ in \eqref{eq:nlsCW0} is equivalent to
changing $\lambda$ in  \eqref{eq:nls},  we see here a behavior
analogous to Test~8 (Figure~\ref{im:supercr}), where similar data are
used. The value $C=1.8$ is very close to (actually slightly below)
the potential energy level  where the blow-up changes from two point
blow-up to one point blow-up, and non-monotonicity can be observed. 
\end{rema}
  
\smallbreak



\section{Numerical test, damped NLS}
\label{sec:numerdamped}

\subsection{Critical power}
We now turn to \eqref{eq:damped}. 
\smallbreak

\noindent{\bf Test in one space dimension}

\begin{figure}
  \center
  \includegraphics[width=0.5\textwidth]{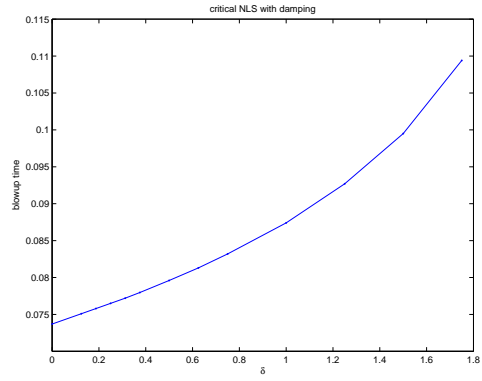}
  \caption{\label{im:dampedsing}
  Blow-up time with varying damping constant $\delta$ (Test~16).
 }
\end{figure}

{\bf Test 16.} 
For  \eqref{eq:damped} in space dimension one, we first use the 
single Gaussian data \eqref{eq:singledata}. 
The scale $C=2$ was chosen so that $\|u_0\|^2_{L^2}=5.013$,
the smallest discretization parameters used are
$Np=2^{12}$ mesh points and $\Delta t = 2.5 \cdot 10^{-6}$.
The blow-up time with respect to changing $\delta$ is shown in
Figure~\ref{im:dampedsing}. 
The blow-up time is decreasing monotonically with $\delta$, in 
accordance with the arguments of \cite{Fibich01}.
For $\delta > 1.75$, blow-up is prevented.

{\bf Test 17.} 
Next we use the ``two-hump''  data 
\eqref{eq:twohpdata}.
The data scale was chosen as $C=5$ so that $\|u_0\|^2_{L^2}=11.969$.

The finest discretization parameters used in
this test are
$Np=2^{14}$ meshpoints and  $\Delta t = 2.5 \cdot 10^{-7}$.

\begin{figure}
  \center
  \includegraphics[width=0.6\textwidth]{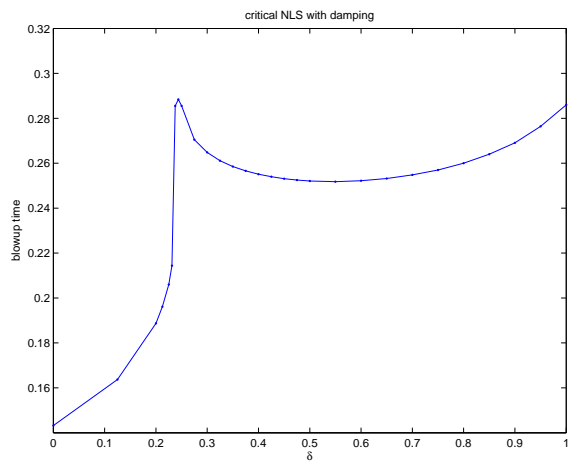}
  \caption{\label{imag_damped}
  Blow-up time with varying damping constant $\delta$ (Test~17).
 }
\end{figure}

The blow-up time with respect to changing $\delta$ is shown in
Figure~\ref{imag_damped}.  
It can be seen that the blow-up time is not monotonously increasing
with $\delta$. The effect is somehow more pronounced than in 
the case of Equation~\eqref{eq:nls}.
\smallbreak

\smallskip
\noindent {\bf Test in two space dimensions}
\smallskip

{\bf Test 18.} For the two-dimensional case of \eqref{eq:damped}, we use the 
data \eqref{eq:tdth_data} with $C=11$, so that the initial 
mass is $\|u_0\|^2_{L^2} = 37.04$. This is above the
minimal value necessary for blow-up~(\cite{Fibich01}).
The finest discretization parameters used in
this test are
$Np=2^{11}$ meshpoints  and  $\Delta t = 5.0 \cdot 10^{-6}$. 
The discretization domain is $[-8,8]^2$.

\begin{figure}
  \center
  \includegraphics[width=0.6\textwidth]{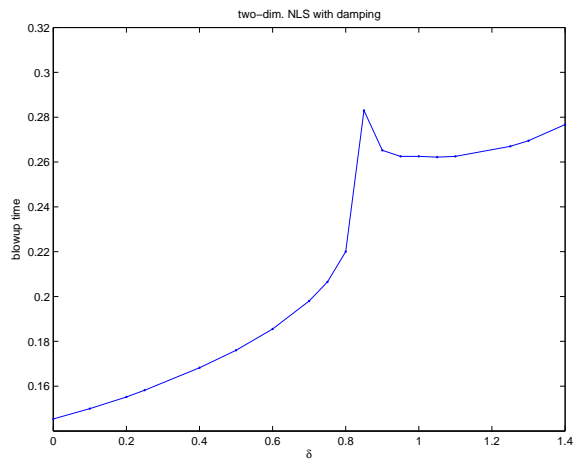}
  \caption{\label{im:twod_damped}
  Blow-up time with varying damping constant $\delta$ (Test~18).
 }
\end{figure}

The blow-up time in dependence of $\delta$ is shown in 
Figure~\ref{im:twod_damped}. For $\delta > 1.4$, blow-up is prevented.

\subsection{Supercritical Power}
$\ $

{\bf Test 19.} We tested equation  \eqref{eq:damped} in space dimension one
with $\sigma = 3$.
The data are chosen as \eqref{eq:twohpdata} with $C=3.8$,
so $\|u_0\|_{L^2}^2 = 3.53$.
The discretization parameters used in
this test are
$Np=2^{12}$ meshpoints and  $\Delta t = 5.0 \cdot 10^{-6}$ or 
 $\Delta t = 10^{-5}$. 
The discretization domain is $[-8,8]$.

\begin{figure}
  \center
  \includegraphics[width=0.6\textwidth]{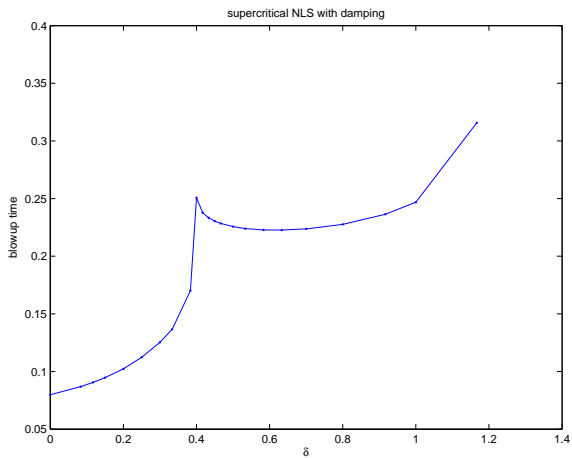}
  \caption{\label{im:supercr_damped}
  Blow-up time with varying damping constant $\delta$,
  supercritical power (Test~19).
 }
\end{figure}

The blow-up times with respect to the damping constant $\delta$
are shown in Figure~\ref{im:supercr_damped}.
Also in this case monotonicity is not true.
For $\delta > 1.17$, blow-up is prevented.
In this case, the occurrence of non-monotonicity is very sensitive to the
size of the data, as already observed in Test 16.
If we choose $C = 3.6$
or $4.0$, instead of $C=3.8$, we observe monotonicity.

\begin{figure}
  \center
  \includegraphics[width=0.5\textwidth]{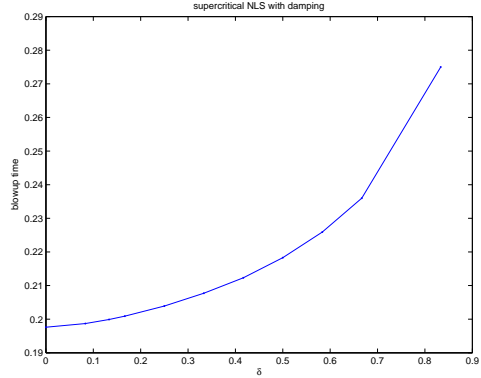}
  \caption{\label{im:supercr_low}
  Blow-up time with varying damping constant $\delta$,
  supercritical power (Test~20).
 }
\end{figure}

{\bf Test 20.} We repeat Test 19 with a different scale of the data,
we take  $C=3.6$ in \eqref{eq:twohpdata},
so that $\|u_0\|_{L^2}^2 = 3.17$.
The discretization parameters used in
this test are
$ Np=2^{15}$ mesh points and  $\Delta t = 5.0 \cdot 10^{-6}$. 
The discretization domain is $[-8,8]$.
The blow-up times with respect to the damping constant $\delta$
are shown in Figure~\ref{im:supercr_low}.
We can see that the blow-up is monotonous in $\delta$. The
blow-up always happens a single point at the origin, and for 
$\delta > 0.83$, blow-up is prevented.

\bigskip

All tests were done with both the TSSP 
and the Relaxation scheme (RS). The results agree, as an example
we showed the comparison of the schemes in figure \ref{im:quadsc_nom}. 
By using two numerical schemes with different discretization
approaches, the possibility of
 observing just numerical defects introduced
by a particular discretization method can be excluded. 
The effects of non-monotonicity found here 
can be seen with two inherently different numerical methods,
which is an indication that the observations are not related to
numerical errors but indeed analytical properties.

\section{Conclusion}
\label{sec:conclusion}

In this paper, we have addressed numerically the question of the
dependence of the blow-up time for solutions to nonlinear
Schr\"odinger equations upon one specific parameter, in three cases:
\begin{itemize}
\item Dependence upon the coupling constant $\l$, for fixed $n,\si$
  and $u_0$:
\begin{equation*}
  i\d_t u +\Delta u = -\lambda |u|^{2\si}u\ , \ (t,x)\in
  \R_+\times \R^n\quad ;\quad
  u_{\mid t =0 } = u_0\, .
\end{equation*}
\item Dependence upon the magnitude of a quadratic oscillation
  introduced in the initial data: only $a$ varies in the following
  equation, 
\begin{equation*}
  i\d_t w + \Delta w = - |w|^{2\si} w  \quad
  ; \quad w_{\mid t=0} = u_0(x)e^{-i\frac{|x|^2}{4a}}\, .
\end{equation*}
\item Dependence upon the strength of the damping $\delta \geq 0$ in 
\begin{equation*}
  i\d_t \psi + \Delta \psi = - |\psi|^{2\si} \psi-i\delta
  \psi \quad ; \quad
  \psi_{\mid t=0} = u_0\, .
\end{equation*}
\end{itemize}
In the $L^2$ super-critical case $\si >\frac{2}{n}$, the tests we
performed highly suggest that in either of the above three cases, the
blow-up time is not monotonous with respect to the variation of the
parameter considered. 
\smallbreak

In the $L^2$ critical case $\si =\frac{2}{n}$, there is apparently no
monotonicity in the first and in the third problem. In the second one
however, our tests agree with the analytical result: there is
monotonicity of the blow-up time with respect to $a$, as recalled in
Proposition~\ref{prop:quadth}.
\smallbreak

We used two numerical methods (time-splitting
spectral method  and a relaxation method), for the results observed to
be more convincing. We may say that they are, since the two methods
yield the same results (and not only just similar results, see
 Figure~\ref{im:quadsc_nom}).  
\smallbreak

Note that in some cases where we observed monotonicity
reversal, the slope of the blow-up time/varying parameter curve may be 
rather steep near the monotonicity breakup. Compare
Figure~\ref{im:tuponedown} near $\lambda=2.7$ with Figure~\ref{im:test1}. In
the quadratic oscillations case, Figure~\ref{im:quadsc_nom} shows a
similar feature near $a=0.7$. For the equation with a damping term,
compare Figures~\ref{imag_damped}  and \ref{im:twod_damped} with
Figure~\ref{im:dampedsing}. This suggests that the underlying
mechanism causing the monotonicity reversal is strongly nonlinear. It
seems inappropriate to speak of instability though (for instance in
Figure~\ref{im:tuponedown}, we saw that there is one point on the
steep line after the local maximum). In these examples, the blow-up
time changes rather fast compared to the modification of  the
$L^2$-norm of the initial data (see also Remark~\ref{rema:C}). Yet, the
dependence of the blow-up time might also involve other analytical
quantities, or  geometrical features.  
\smallbreak

All the numerical counter-examples to monotonicity that
we found contain a somehow nontrivial profile, inasmuch as the initial
datum is formed of two (or more) humps. The question of an analytical
justification remains open and challenging in these cases.


\appendix
\section{Proofs}
\label{sec:proofs}

\subsection{Proof of Proposition~\ref{prop:local}}
\label{sec:proplocal}

It is based on a fixed
point argument using Strichartz estimates, which we first recall.
\begin{defin}\label{def:adm}
 A pair $(q,r)$ is {\bf admissible} if $2\leq r
  \leq\frac{2n}{n-2}$ (resp. $2\leq r\leq \infty$ if $n=1$, $2\leq r<
  \infty$ if $n=2$)  
  and 
$$\frac{2}{q}=\delta(r):= n\left( \frac{1}{2}-\frac{1}{r}\right).$$
\end{defin}
\begin{lem}[Strichartz estimates, see
  e.g. \cite{Yajima87,GV92,KT}]\label{lem:strichartz} 
Let $(q,r)$, $(q_1,r_1)$ and~$ (q_2,r_2)$ be admissible pairs. Denote
$U(t) := e^{it\Delta}$.\\
$1.$ There exists $C_r$ such that for any $\varphi \in L^2(\R^n)$:
\begin{equation}\label{eq:strich}
    \left\| U(.)\varphi \right\|_{L^q(\R;L^r)}\leq C_r \|\varphi
    \|_{L^2}\, .
  \end{equation}
$2.$ There exists $C_{r_1,r_2}$ such that for any interval 
    $I$ and any $F\in L^{q'_2}(I;L^{r'_2})$:
\begin{equation}\label{eq:strichnl}
      \left\| \int_{I\cap\{s\leq
      t\}} U(t-s)F(s)ds 
      \right\|_{L^{q_1}(I;L^{r_1})}\leq C_{r_1,r_2} \left\|
      F\right\|_{L^{q'_2}(I;L^{r'_2})}\, .
    \end{equation}
\end{lem}
The proof of Proposition~\ref{prop:local} relies on the following
straightforward lemma:
\begin{lem}\label{lem:alg}
Let $r=s =2\si +2$, and $q=\frac{4\si+4}{n\si}$, so that the pair $(q,r)$ is
admissible. Define $k$ by
$$k=\frac{2\si (2\si +2)}{2-(n-2)\si}.$$
Then $k$ is finite, and the following algebraic identities hold:
\begin{equation*}
\frac{1}{r'}= \frac{1}{r}+\frac{2\si}{s}\quad ;\quad 
\frac{1}{q'}= \frac{1}{q}+\frac{2\si}{k}\ \cdot
\end{equation*}
\end{lem}
\begin{proof}[Proof of Proposition~\ref{prop:local}]
For $A\in \{ {\rm Id},\nabla_x\}$, denote $R_A := \| A
u_0\|_{L^2}$. With the notations of Lemmas~\ref{lem:strichartz} and 
\ref{lem:alg}, define:
\begin{equation*}
\begin{aligned}
 X_T := \big\{ & u \in C([0,T[;H^1)\ ; \ Au\in
 L^\infty([0,T[;L^2)\cap L^q([0,T[;L^r), \\ 
 \|Au&\|_{L^\infty([0,T[;L^2)}\leq 
2R_A\ ,\ 
 \|Au\|_{L^q([0,T[;L^r)}\leq 2C_{2\si}R_A, \ 
\forall A\in \{ {\rm Id},\nabla_x\} \big\} \, .
\end{aligned}
\end{equation*}
Duhamel's formula for \eqref{eq:nls} writes:
\begin{equation}\label{eq:Duhamel}
u(t)=U(t)u_0 +i \lambda \int_0^t
U(t-s)\left(|u|^{2\si}u\right)(s)ds\, .
\end{equation} 
Denote $F(u)$ the right hand side of \eqref{eq:Duhamel}. We prove that
$F$ maps $X_T$ into itself and is a contraction (in a weaker metric)
provided that \eqref{eq:cond} is satisfied, for $\e>0$ sufficiently
small. Denote $L^a_T :=L^a([0,T[)$. From Lemmas~\ref{lem:strichartz} and 
\ref{lem:alg}, we have, for $A\in \{ {\rm Id},\nabla_x\}$:
\begin{equation}\label{eq:est1}
\begin{aligned}
  \left\|A\(F(u)\) \right\|_{L^\infty_T(L^2)}& \leq
  \left\|Au_0 \right\|_{L^2} +C \lambda \left\|A\(|u|^{2\si}u\)
  \right\|_{L^{q'}_T(L^{r'})}\\
 &\leq
  \left\|Au_0 \right\|_{L^2} +C\lambda  \left\|u
  \right\|^{2\si}_{L^{k}_T(L^{s})}\left\|Au
  \right\|_{L^{q}_T(L^{r})}\,,
\end{aligned}
\end{equation}
for some constant $C$ depending only on $n$ and $\si$.
From Gagliardo--Nirenberg inequality, we have:
\begin{equation*}
 \left\|u(t)
  \right\|_{L^{s}}\leq C(n,s) \left\|u(t)
  \right\|_{L^{2}}^{1-\delta(s)}\left\|\nabla_xu(t)
  \right\|_{L^{2}}^{\delta(s)} \, , 
\end{equation*}
where $\delta(s)$ is given in Definition~\ref{def:adm}. We infer, for
$u\in X_T$:
\begin{equation*}
 \left\|A\(F(u)\) \right\|_{L^\infty_T(L^2)} \leq R_A + C\lambda
 T^{2\si /k}
 R_{{\rm Id}}^{2\si (1-\delta(s))} R_{\nabla_x}^{2\si \delta(s)} R_A
 \, ,
\end{equation*}
where from now on, we denote by $C$ all the constants which do not
depend on relevant parameters. 
Using Lemma~\ref{lem:strichartz}, we have similarly:
\begin{equation*}
 \left\|A\(F(u)\) \right\|_{L^q_T(L^r)} \leq C_rR_A + C\lambda T^{2\si
 /k}
 R_{{\rm Id}}^{2\si (1-\delta(s))} R_{\nabla_x}^{2\si \delta(s)} R_A
 \, ,
\end{equation*}
for some other constant $C$ depending only on $n$ and $\si$. We
therefore have stability for
\begin{equation}\label{eq:cond0}
  \lambda T^{2\si/k}
 R_{{\rm Id}}^{2\si (1-\delta(s))} R_{\nabla_x}^{2\si \delta(s)}\ll
 1\, .
\end{equation}
We have contraction in the weaker metric $L^q([0,T[;L^r)$ under a
similar condition:
\begin{equation*}
\begin{aligned}
    \big\| F(u_2)-F(u_1)
    \big\|_{L^q_T(L^r)}  & \leq
C_{r,r}\lambda\left\| \left(|u_2|^{2\si} u_2 - |u_1|^{2\si}
    u_1\right) \right\|_{L^{q'}_T(L^{r'})}\\
&\leq C\lambda\left(
    \|u_1\|^{2\si}_{L^k_T(L^s)}+
\|u_2\|^{2\si}_{L^k_T(L^s)}\right) 
\|u_2-u_1\|_{L^q_T(L^r)}\\
&\leq C\lambda T^{2\si/k}
 R_{{\rm Id}}^{2\si (1-\delta(s))} R_{\nabla_x}^{2\si \delta(s)}
\|u_2-u_1\|_{L^q_T(L^r)}.
\end{aligned}
\end{equation*}
If \eqref{eq:cond0} is satisfied, then $F$ has a unique
fixed point in $X_T$. Since  \eqref{eq:cond0} is nothing else
but \eqref{eq:cond}, this proves the first part of
Proposition~\ref{prop:local}. 

When $u_0\in \Sigma$, one can prove that the $\Sigma$-regularity is
conserved along the time evolution by now considering $A\in \{ {\rm
  Id},\nabla_x,J(t)\}$, where $J(t)=x+2it \nabla_x$ is the Galilean
operator. (This operator commutes with the group $U(t)$ and acts on the
nonlinearity we consider like a derivative.) Finally, we refer to
\cite{Caz} for the conservation laws.
\end{proof}

\subsection{Proof of Corollary~\ref{cor:improved}}
\label{sec:corimproved}

In the above proof of local existence, we did not use the conservation
laws. We now take them into account.

  Consider Duhamel's formula \eqref{eq:Duhamel}, and seek an
  $L^\infty_t(L^2_x)$-bound  for $\nabla_x u$. For $A\in \{ {\rm
  Id},\nabla_x\}$, we have like above:
\begin{align*}
  \left\|Au \right\|_{L^q_T(L^r)}& \leq
  C_r \left\|Au_0 \right\|_{L^2} +C\lambda  \left\|u
  \right\|^{2\si}_{L^{k}_T(L^{s})}\left\|Au
  \right\|_{L^{q}_T(L^{r})}\,.
\end{align*}
Recall that $s=2\si +2$; from the conservation of energy,
\begin{equation*}
  \left\|u(t)
  \right\|_{L^{s}}\leq C \l^{-1/s}\(\| \nabla_x u(t)\|_{L^2}^{2/s}+
  |E| \)\leq C\(\l^{-1/s}\| \nabla_x u(t)\|_{L^2}^{2/s}+ 1\) , 
\end{equation*}
for $\l\geq 1$, where $C$ depends on $u_0$ and $\si$, but not on
$\l$. 
We deduce, for any $t\in [0,T]$:
\begin{equation}\label{eq:abs}
  \begin{aligned}
    \left\|Au \right\|_{L^q_t(L^r)} \leq & \ 
  C_r \left\|Au_0 \right\|_{L^2}\\
&+CT^{2\si /k}\(\lambda^{\frac{\si}{\si +1}}\| \nabla_x
  u\|_{L^\infty_t(L^2)}^{4\si/s}+ 1 \)\left\|Au
  \right\|_{L^{q}_t(L^{r})}.
  \end{aligned}
\end{equation}
Now assume that \eqref{eq:cond2} is satisfied for some $\e$ to be
fixed later. Then by continuity (see
Proposition~\ref{prop:local}), there exists $\tau>0$ such that
\begin{equation}
  \label{eq:solong}
   T^{2-(n-2)\si} \(
\lambda^{2\si}\|\nabla_x u\|_{L^\infty_t (L^2)}^{4\si}+1\)\leq
\(2^{4\si}+1\)\e\, , 
\end{equation}
for $0\leq t\leq \tau$.
Choosing $0<\e\leq \e_0(n,\si)$, the last term of \eqref{eq:abs} can be
absorbed by the left hand side for $0\leq t\leq \tau$:
\begin{equation*}
  \left\|Au \right\|_{L^q_t(L^r)} \leq
  2 C_r \left\|Au_0 \right\|_{L^2}.
\end{equation*}
Strichartz inequalities and the above estimate yield, for $t\in [0,T]$
such that \eqref{eq:solong} holds:
\begin{equation*}
  \begin{aligned}
  \left\|Au \right\|_{L^\infty_t(L^2)} &\leq
   \left\|Au_0 \right\|_{L^2}
  +CT^{2\si /k}\(\lambda^{\frac{\si}{\si +1}}\| \nabla_x
  u\|_{L^\infty_t(L^2)}^{4\si/s}+ 1 \)\left\|Au
  \right\|_{L^{q}_t(L^{r})}\\
&\leq
   \left\|Au_0 \right\|_{L^2}
  + C'T^{2\si /k}\(\lambda^{\frac{\si}{\si +1}}\| \nabla_x
  u\|_{L^\infty_t(L^2)}^{4\si/s}+ 1 \)\left\|Au_0 \right\|_{L^2}.
  \end{aligned}
\end{equation*}
By the conservation of mass, this estimate is interesting
only when $A=\nabla_x$. 
Up to choosing $\e$ even smaller, we see that
as long as \eqref{eq:solong} is satisfied:
\begin{equation*}
  \left\|\nabla_x u \right\|_{L^\infty_t(L^2)} \leq
   2\left\|\nabla_x u_0 \right\|_{L^2}.
\end{equation*}
Now suppose that \eqref{eq:solong} is not satisfied for $t=T$: there
exists a minimal time $t^*<T$ such that
\begin{align*}
   T^{2-(n-2)\si} \(\lambda^{2\si}
\|\nabla_x u\|_{L^\infty_{t^*} (L^2)}^{4\si}+1\)&= \(2^{4\si}+1\)\e\ ,\\
\  \left\|\nabla_x u \right\|_{L^\infty_{t^*}(L^2)}& \leq
   2\left\|\nabla_x u_0 \right\|_{L^2}.
\end{align*}
The latest inequality yields, along with \eqref{eq:cond2}:
\begin{equation*}
   T^{2-(n-2)\si} \(\lambda^{2\si}
\|\nabla_x u\|_{L^\infty_{t^*} (L^2)}^{4\si}+1\)\leq 
T^{2-(n-2)\si} \( \lambda^{2\si}2^{4\si}
\|\nabla_x u_0\|_{L^2}^{4\si}+1\)\leq 2^{4\si} \e.
\end{equation*}
This contradicts the definition of $t^*$. Therefore, \eqref{eq:solong}
is  satisfied for $t=T$, which completes the proof of
Corollary~\ref{cor:improved}.

\subsection{Proof of Proposition~\ref{prop:quadth}}
\label{sec:propquadth}
For $\alpha \in \R$, define:
\begin{equation*}
J_\alpha(t) = x + 2i(t-\alpha)\nabla_x=
2i(t-\alpha)e^{i\frac{|x|^2}{4(t-\alpha)}}
\nabla_x\left(e^{-i\frac{|x|^2}{4(t-\alpha)}} \cdot \right) \, .
\end{equation*}
Denote
\begin{equation*}
\phi(t)=h(t)^{n\si -2}= \left(1-\frac{t}{a} \right)^{n\si -2}\, .
\end{equation*}
Besides the conservations of mass for $u$ and $v$, and of energy for
$u$, we have the following evolution laws:
\begin{align*}
  \frac{d}{dt}&\Big(\|J_\alpha(t) u\|_{L^2}^2 -4\frac{(t-\alpha)^2}{\si
+1}
\| u(t)\|_{L^{2\si+2}}^{2\si+2}\Big) 
=4\frac{n\si-2}{\si
+1}(t-\alpha)\| u(t)\|_{L^{2\si+2}}^{2\si+2};  \\
\frac{d}{dt}&\Big(\|\nabla_x v(t)\|_{L^2}^2 -\frac{\phi(t)}{\si +1}
\| v(t)\|_{L^{2\si+2}}^{2\si+2}\Big) = \frac{-\phi
'(t)}{\si +1}\| v(t)\|_{L^{2\si+2}}^{2\si+2};\\ 
\frac{d}{dt}&\left(\|J_\alpha(t) v\|_{L^2}^2
-4\frac{(t-\alpha)^2\phi(t)}{\si +1}
\| v(t)\|_{L^{2\si+2}}^{2\si+2}\right) =\\
&=4\frac{n\si-2}{\si
+1}(t-\alpha)\phi(t)\| v(t)\|_{L^{2\si+2}}^{2\si+2}
-\frac{4}{\si
+1}(t-\alpha)^2\phi'(t)
\| v(t)\|_{L^{2\si+2}}^{2\si+2} \, .
\end{align*}

\begin{lem}\label{lem:blow0}
Let $u_0\in\Sigma$, $2/n\leq \si<2/(n-2)$ and $T>0$. The following
assertions are equivalent:
\begin{itemize}
\item[$(1)$] $u$ blows up at time $T$.
\item[$(2)$] $\|\nabla_x u(t)\|_{L^2} \to +\infty$ as $t\to T$.
\item[$(3)$] $\| u(t)\|_{L^{2\si+2}} \to +\infty$ as $t\to T$.
\item[$(4)$] For any $\alpha \not = T$, $\| J_\alpha (t)u\|_{L^2}
\to +\infty$ as $t\to T$.
\end{itemize}
\end{lem}
\begin{proof}[Proof of Lemma~\ref{lem:blow0}]
The equivalence between the first two points follows from
Proposition~\ref{prop:local}. The second point is equivalent to the
third one by conservation of the energy.

For  $ \alpha \not = T$, Gagliardo--Nirenberg inequalities and the
second writing for $J_\alpha$ yield:
\begin{equation*}
\|u(t)\|_{L^{2\si+2}}^{2\si+2} \leq  \frac{C(n,\si)}{|t-\alpha|^{n\si}}\|
u(t)\|_{L^2}^{2-(n-2)\si} \| J_\alpha (t)u\|_{L^2}^{n\si}\, .
\end{equation*}
We infer (3)$\Rightarrow$(4) by the conservation of mass. 
Finally, if $\| u(t)\|_{L^{2\si+2}}$ remains bounded as $t\to
T$, then integrating the first evolution law above 
shows that $\|J_\alpha(t) u\|_{L^2}$ remains bounded as $t\to
T$. 
\end{proof}
A similar result holds for $v$. Notice that if $a>0$, then for
$0<t<a$, $\phi$  is smooth (bounded) and does not cancel. If $a<0$,
the same holds in a compact of  $\R_+$. 
\begin{lem}\label{lem:explv}
Let $u_0\in\Sigma$, $2/n\leq \si<2/(n-2)$ and $T>0$, with
$T<a$ if $a>0$. The following
assertions are equivalent:
\begin{itemize}
\item[$(1)$] $v$ blows up at time $T$.
\item[$(2)$] $\|\nabla_x v(t)\|_{L^2} \to +\infty$ as $t\to T$.
\item[$(3)$] $\| v(t)\|_{L^{2\si+2}} \to +\infty$ as $t\to T$.
\item[$(4)$] For any $\alpha \not = T$, $\| J_\alpha (t)v\|_{L^2}
\to +\infty$ as $t\to T$.
\end{itemize}
\end{lem}
\begin{proof}
The equivalence between (1) and (2) follows like for
Proposition~\ref{prop:local}, since $\phi$ is bounded from
above and from below by positive constants.

The other equivalences follow like above, thanks to
Gagliardo--Nirenberg inequalities.   
\end{proof}
Then Proposition~\ref{prop:quadth} follows from Lemma~\ref{lem:explv}
and the identity:
\begin{equation*}
\| J_a(t)v\|_{L^2} =2\left\|
\nabla_x u\left(\frac{at}{a-t}\right)\right\|_{L^2}\, .
\end{equation*}

\subsection{Proof of Proposition~\ref{prop:damped}}
\label{sec:propdamped}
It consists of a slight adaptation of the proof of
Proposition~\ref{prop:local}. Considering \eqref{eq:damped}, the
analog of \eqref{eq:est1} yields:
\begin{equation*}
\begin{aligned}
  \left\|A\(F(u)\) \right\|_{L^\infty_T(L^2)}& \leq
  \left\|Au_0 \right\|_{L^2} +C \lambda \left\|e^{-2\si \delta
  t}A\(|u|^{2\si}u\) 
  \right\|_{L^{q'}_T(L^{r'})}\\
 &\leq
  \left\|Au_0 \right\|_{L^2} +C\lambda  \left\|e^{-\delta
  t} u
  \right\|^{2\si}_{L^{k}_T(L^{s})}\left\|Au
  \right\|_{L^{q}_T(L^{r})}\,,
\end{aligned}
\end{equation*}
where the notation $L^p_T$ now stands for $L^p([0,T[)$. Using
Gagliardo--Nirenberg inequality, we have:
\begin{equation*}
  \left\|e^{-\delta
  t} u
  \right\|_{L^{k}_T(L^{s})} \leq C \left\|e^{- \delta
  t} 
  \right\|_{L^{k}_T} \left\|u
  \right\|_{L^{\infty}_T(L^2)}^{1-\delta(s)} \left\|\nabla_x u
  \right\|_{L^{\infty}_T(L^2)}^{\delta(s)}  .
\end{equation*}
Mimicking the proof of Proposition~\ref{prop:local}, we can apply a
fixed point argument provided that:
\begin{equation*}
  \left\|e^{-\delta
  t} 
  \right\|_{L^{k}_T}
 R_{{\rm Id}}^{1-\delta(s)} R_{\nabla_x}^{ \delta(s)}\ll
 1\, . 
\end{equation*}
This property is a consequence of the stronger one:
\begin{equation*}
  R_{{\rm Id}}^{1-\delta(s)}
  R_{\nabla_x}^{ \delta(s)}\ll \delta  ,
\end{equation*}
which is the condition stated in
Proposition~\ref{prop:damped}, whose proof is now complete.
\begin{rema*}
  This proof yields similar results for a larger class of damping. For
  $\delta >0$ and $F\in L^1(\R_+)$ a positive function, define:
  \begin{equation*}
    H(t,\delta) = -\delta \frac{F'(\delta t)}{F(\delta t)}\, .
  \end{equation*}
Then for $\delta$ satisfying \eqref{eq:conddamp} with another constant
$C$ (depending of $F$), the solution of
\begin{equation*}
  i\d_t \psi + \Delta \psi = - |\psi|^{2\si} \psi-iH(t,\delta) 
  \psi \, ,\ x\in 
  \R^n \quad ; \quad
  \psi_{\mid t=0} = u_0\, ,
\end{equation*}
does not blow up in the future. In Prop.~\ref{prop:damped}, we
considered the case $F(y)=e^{-y}$. Considering $F(y)=(1+y)^{-\alpha}$
with $\alpha >1$ yields a damping function
\begin{equation*}
    H(t,\delta) = \alpha \frac{\delta}{1+\delta t}\, , 
  \end{equation*}
which goes to zero as $t \to +\infty$. Even though the damping is less
and less strong, \eqref{eq:conddamp} ensures global existence in $H^1$
in the future.
\end{rema*}


\providecommand{\href}[2]{#2}

\end{document}